\documentclass[11pt,twoside]{article}

\usepackage{latexsym}
\usepackage{ifthen}
\usepackage{cite}

\input amssym.def
\input amssym.tex

\newcommand{\In}{\infty}
\newcommand{\de}{\delta}
\newcommand{\eps}{\epsilon}

\newcommand{\sig}{\sigma}
\newcommand{\om}{\omega}
\newcommand{\De}{\Delta}
\newcommand{\La}{\Lambda}
\newcommand{\RR}{{\Bbb R}}

\newcommand{\ZZ}{{\Bbb Z}}
\newcommand{\sX}{{\cal X}}
\newcommand{\sO}{{\cal O}}
\newcommand{\sF}{{\cal F}}
\newcommand{\sC}{{\cal C}}
\newcommand{\sK}{{\cal K}}
\newcommand{\sI}{{\cal I}}
\newcommand{\sA}{{\cal A}}
\newcommand{\sT}{{\cal T}}
\newcommand{\sB}{{\cal B}}
\newcommand{\sP}{{\cal P}}
\newcommand{\be}{{\bf e}}
\newcommand{\bff}{{\bf f}}
\newcommand{\bg}{{\bf g}}
\newcommand{\bm}{{\bf m}}
\newcommand{\bp}{{\bf p}}
\newcommand{\bt}{{\bf t}}
\newcommand{\bv}{{\bf v}}
\newcommand{\bw}{{\bf w}}
\newcommand{\bx}{{\bf x}}
\newcommand{\by}{{\bf y}}
\newcommand{\bP}{{\bf P}}

\newcommand{\rank}{{\rm rank}}
\newcommand{\vol}{{\rm vol}}
\newcommand{\supp}{{\rm supp}}

\newcommand{\beql}[1]{\begin{equation}\label{#1}}
\newcommand{\eeq}{\end{equation}}

\catcode`\@=11
\renewcommand{\section}{\setcounter{equation}{0}%
\@startsection{section}{1}{0ex}{-4.3ex}{1pt}{\it}}
\renewcommand{\paragraph}%
{\@startsection{paragraph}{5}{0ex}{1ex}{-1em}{\it}}
\renewcommand{\@biblabel}[1]{[\bf #1]}
\renewcommand{\@cite}[2]{%
[{\bf #1}\ifthenelse{\boolean{@tempswa}}{,#2}{}]}

\renewcommand{\@makefnmark}{\mbox{$\@thefnmark$~}}
\renewcommand{\@makefntext}[1]{\noindent\@makefnmark #1}
\newcounter{fn}[page]
\newcommand{\fn}[1]{\addtocounter{fn}{2}\footnote{#1}\addtocounter{fn}{-1}} 

\catcode`@=12

\makeatletter
\def\eqalignno#1{\displ@y \ta {\bf s} kip\@centering
  \halign to\displaywidth{\hfil$\@lign\displaystyle{##}$\ta {\bf s} kip\z@skip
    & $\@lign\displaystyle{{}##}$\hfil\ta {\bf s} kip\@centering
    & \llap{$\@lign##$}\ta {\bf s} kip\z@skip\crcr
    #1\crcr}}
\makeatother

\makeatletter
\def\@sect#1#2#3#4#5#6[#7]#8{\ifnum #2>\c@secnumdepth
     \def\@svsec{}\else 
     \refstepcounter{#1}\edef\@svsec{\csname the#1\endcsname.\hskip .75em }\fi
     \@tempskipa #5\relax
      \ifdim \@tempskipa>\z@ 
	\begingroup #6\relax
	  \@hangfrom{\hskip #3\relax\@svsec}{\interlinepenalty \@M #8\par}%
	\endgroup
       \csname #1mark\endcsname{#7}\addcontentsline
	 {toc}{#1}{\ifnum #2>\c@secnumdepth \else
		      \protect\numberline{\csname the#1\endcsname}\fi
		    #7}\else
	\def\@svsechd{#6\hskip #3\@svsec #8\csname #1mark\endcsname
		      {#7}\addcontentsline
			   {toc}{#1}{\ifnum #2>\c@secnumdepth \else
			     \protect\numberline{\csname the#1\endcsname}\fi
		       #7}}\fi
     \@xsect{#5}}
\def\@theorem#1#2{\it \trivlist \item[\hskip \labelsep{\bf #1\ #2.}]}
\makeatother  

\newtheorem{Theorem}{\sc Theorem}[section]
\newtheorem{Corollary}{\sc Corollary}[section]
\newtheorem{Lemma}{\sc Lemma}[section]
\newtheorem{Conjecture}{\sc Conjecture}[section] 
\newtheorem{Example}{\sc Example}[section] 

\newenvironment{enum}{\begin{enumerate}\setlength{\itemsep}{0pt}%
\setlength{\parsep}{0pt}\vspace{-2mm}}{\vspace*{-2mm}\end{enumerate}}
\newenvironment{enumalph}
{%
\begin{enum}}{\end{enum}}
\newenvironment{enumarabic}
{%
\begin{enum}}{\end{enum}}
\newenvironment{enumroman}
{%
\begin{enum}}{\end{enum}}

\begin{document}

\vspace*{2cm}

\begin{center}
\addtocounter{fn}{2} 
{\LARGE\bf Repetitive Delone sets and quasicrystals}\\ \bigskip
JEFFREY C.~LAGARIAS\dag\ and PETER A.~B.~PLEASANTS\ddag%
\fn{Current address: Department of Mathematics,
The University of Queensland, Queensland 4072, Australia.}
\\ \smallskip
\dag{\em AT\&T Labs---Research, Florham Park, NJ 07932, USA \\
(e-mail: jcl@research.att.com)}\\
\ddag{\em Department of Mathematics and Computing Science,
University of the South Pacific,\\
Suva, Fiji \\
(e-mail: pabp@maths.uq.edu.au)} 
\vspace{2cm}
\end{center}

\paragraph{Abstract.}  This paper studies the problem of characterizing
the simplest aperiodic discrete point sets, using invariants based on
topological dynamics.
A Delone set of finite type is a Delone set $X$ such that $X-X$ is
locally finite.
Such sets are characterized by their patch-counting function $N_X(T)$
of radius $T$ being finite for all $T$.
We formulate conjectures relating slow growth of the patch-counting
function $N_X(T)$ to the set $X$ having a non-trivial translation symmetry.

A Delone set $X$ of finite type is repetitive if there is a function
$M_X(T)$ such that every closed ball of radius $M_X(T) + T$ contains
a complete copy of each kind of patch of radius $T$ that occurs in $X$.
This is equivalent to the minimality of an associated topological
dynamical system  with $\RR^n$-action.
There is a lower bound for $M_X(T)$ in terms of $N_X(T)$, namely
$M_X(T) \ge c(N_X (T))^{1/n}$ for some positive constant $c$ depending
on the Delone set constants $r,R$, but there is no general upper bound
for $M_X(T)$ purely in terms of $N_X(T)$.
The complexity of a repetitive Delone set $X$ is measured
by the growth rate of its repetitivity function $M_X (T)$.
For example, the function $M_X (T)$ is bounded if and only
if $X$ is a periodic crystal.
A set $X$ is {\em linearly repetitive} if $M_X (T) = O(T)$
as $T \to \In$ and is {\em densely repetitive} if
$M_X(T) = O(N_X(T))^{1/n}$  as $T \to \In$.
We show that linearly repetitive sets and densely repetitive
sets have strict uniform patch frequencies, i.e.\ the associated
topological dynamical system is strictly ergodic.
It follows that such sets are diffractive, in the sense of having
a well-defined diffraction measure.
In the reverse direction, we construct a repetitive Delone set 
$X$ in $\RR^n$ which has 
$M_X (T) = O(T (\log T)^{2/n} (\log\log\log T )^{4/n})$,
but does not have uniform patch frequencies.
Aperiodic linearly repetitive sets have many claims to be 
the simplest class of aperiodic sets and we propose considering
them as a notion of `perfectly ordered quasicrystals'.

\section{Introduction}

The most ordered discrete sets in $\RR^n$ are those with crystalline order, 
which have a full rank lattice of global translation symmetries. Such sets
in $\RR^3$ provide good models for the atomic structure of crystals.
Quasicrystalline materials are materials whose x-ray diffraction spectra
have sharp spots, indicating some long-range order in the atomic structure,
but which lack a lattice of periods---typically they exhibit symmetries
which are impossible for  any fully periodic arrangement of the atoms.
A few quasicrystalline materials are thermodynamically stable with apparent
zero entropy and can be grown as `perfectly ordered quasicrystals',
whose long-range atomic order seems to be as good as that of crystals,
even though this atomic order is apparently aperiodic.  One such example is
the Al--Pd--Mn icosahedral qusicrystalline phase, described in \cite{Getal}.
Most quasicrystalline materials are less ordered, and are thought to have
positive entropy and be entropically stabilized.  The atomic structure of
these materials is not known precisely, although there is now extensive
theory and there are various proposed models which give reasonable
agreement with data, see \cite{Jan92,Moo97,Sen95}.  The less-ordered
quasicrystalline materials are modeled with random tilings,
see \cite{He91,RHB98}. 

The general problem motivating this paper is that of characterizing the
`simplest' locally finite point sets $X$ in $\RR^n$ that are aperiodic.
A {\em locally finite set} is a set with only finitely many points in
any bounded region, and an {\em aperiodic}%
\fn{Some authors call such sets {\em non-periodic}.}
set is one with no translation symmetries.
There is an enormous variety of aperiodic sets which have some 
long-range order, see `The Aperiodic Zoo' in Senechal \cite{Sen95}.
Does there exist a nice characterization of the `simplest' aperiodic sets,
which might be used to define a notion of `perfectly ordered quasicrystal'?
We study this question for Delone sets and our notions of simplicity
involve constraints on the patches in the set and invariants from 
topological dynamics. 

Our definitions are inclusive and allow sets with translation symmetries.
In other words, we view a fully periodic set, or `crystal', as a  
special kind of `quasicrystal'. 

\paragraph{Definition 1.1.}
A {\em Delone set}, or $(r,R)${\em-set},
is a discrete set $X$ in $\RR^n$ that has the following two properties.
\begin{enumroman}
\item{\em Uniform discreteness.}
There is an $r > 0$ such that each open ball of radius
$r$ in $\RR^n$ contains at most one point of $X$.
\item{\em Relative density.}
There is an $R > 0$ such that each closed ball of radius
$R$ in $\RR^n$ contains at least one point of $X$.
\end{enumroman}

The classical measure of long-range order of a Delone set 
is its set of global translation symmetries or, more generally,
its set of global symmetries under Euclidean motions.

\paragraph{Definition 1.2.}
The {\em period lattice} of  a Delone set $X\/$ in $\RR^n$ 
is the lattice of translation symmetries  $\La_X$ of $X$,
given by
\beql{eq2003}
\La_X : = \{ \bt \in \RR^n : X + \bt = X \}.
\eeq
It is a free Abelian group, with rank satisfying 
$0 \leq \rank ( \La_X ) \leq n$.

The most regular Delone sets are those with a full set
of translation symmetries, i.e.\ $\rank (\La_X) = n$.
We call such an $X$ an {\em ideal crystal}.  
As noted earlier, a Delone set $X$ is aperiodic if
$\rank (\La_X) = 0$, i.e.\ $\La_X = \{ {\bf 0}\}$. 

In this paper we consider sets that need not have any long-range
order in this classical sense, but which do have a weaker sort
of order described by restrictions on their `local' structure.
In \cite{Lag,Lag98} the first author introduced several classes
of Delone sets $X$ determined by local restrictions on their set
of interpoint vectors $X - X$. 

\paragraph{Definition 1.3.}
(i) A Delone set $X$ in $\RR^n$ is {\em finitely generated\/} if the
additive group $[X - X]$ in $\RR^n$ generated by the set $X - X$ of
interpoint vectors is finitely generated.

(ii) A Delone set $X$ is a {\em Delone set of finite type\/} if
$X - X$ is locally finite, i.e.\ the intersection of $X - X$ with
any bounded set is finite.

(iii) A Delone set $X$ is a {\em Meyer set} if $X - X$ is a
Delone set.
\bigskip

These concepts form a hierarchy, because any Delone set of finite type 
is finitely generated, and any Meyer set is a Delone set of finite type.
Almost all sets studied in quasicrystallography are Delone sets of
finite type (see \cite{Lag}) and this is the class of sets 
that we study in this paper.
The narrower class of Meyer sets was introduced earlier as giving 
useful models for quasicrystalline structures,
see \cite{Mey95,Moo95,Moo97,MooPat95}. 

We consider measures of complexity of a Delone set  
$X$ in $\RR^n$ based on its `local' structure.
The {\em T-patch} of $X$ centered at a point $\bx \in X$ is
\beql{eq101}
\sP_X (\bx ; T) : = X \cap B( \bx ; T),
\eeq
in which $B( \bx ; T)$ denotes the closed ball of a radius $T$
centered at $\bx$.
The {\em T-0atlas} $\sA_X (T)$ of $X$ consists of all $T$-patches
of $X$ translated to the origin:
\beql{eq102}
\sA_X (T) : = \{\sP_X (\bx; T) - \bx : \bx \in X \}.
\eeq
The {\em atlas} $\sA_X $ of $X$ is the union of all $T$-atlases for $T > 0$.
The structure of the atlas gives the following complexity measure.

\paragraph{Definition 1.4.}
The {\em patch-counting function}, $N_X (T)$, counts the number of
$T$-patches in $X$ up to translation-equivalence, i.e.
$$
N_X(T):=|\sA_X(T)|.
$$

For most Delone sets the patch-counting function is infinite for large $T$.
The Delone sets whose patch-counting functions are finite for all
values of $T$ are exactly the Delone sets of finite type and for
such sets the growth rate of the patch-counting function $N_X (T)$ 
provides a quantitative measure of the complexity of $X$.  

The patch-counting function can be viewed as a generalization of the
word-counting function, or {\em permutation number} $p(T)$, that is studied
in the combinatorics of words, and also generalizes the {\em (symbolic)
complexity function} studied in symbolic dynamics, as in Ferenczi~\cite{Fe96}.
We give a more detailed overview of these areas at the end of this section.

In \cite{Lag} it was  shown that $N_X(T)$ satisfies the upper bound
\beql{eq105}
N_X (T) \leq \exp (c_0 T^n ) \quad \mbox{for all $T > 0$},
\eeq
with a positive constant $c_0 = c_0(X)$.  One can further show that any
Delone set of finite type has a finite {\em configurational entropy}
\beql{eq105b}
 H_c (X) := \lim_{T \to \In} \frac{1}{T^n} \log N_X(T).
\eeq
(The existence of the limit can be established by a subadditivity
argument.)
At the other extreme, the slowest possible growth rate for $N_X(T)$ is
to be eventually constant, which occurs when $X$ is an ideal crystal.
In this case, the eventual value of $N_X(T)$ is equal to the number of points
of $X\/$ in a fundamental domain for the lattice of translations of $X$.
It appears that sufficiently slow growth of the patch-counting function
$N_X(T)$ forces a set $X$ to have global translation symmetries. In
this direction we formulate the following conjecture.

\begin{Conjecture}\label{Cj11}
Any aperiodic Delone set $X$ in $\RR^n$ satisfies
\beql{200b}
\limsup_{T \to \In} \frac {N_X(T)}{T^n} > 0.
\eeq
\end{Conjecture}

In \S2 we formulate stonger conjectures, and give examples
showing that one cannot replace $\limsup$ by $\liminf$ in
this conjecture, when the dimension $n \geq 3.$

A weak type of translational order for a Delone set is given by the notion
of repetitivity and carries with it an associated complexity measure.
\paragraph{Definition 1.5.}
A Delone set $X$ is {\em repetitive}%
\fn{There is a parallel concept for tilings, where various terms
are used: `almost periodic' in \cite{Rob9394,Sol97}; `tilings with local
isomorphism' in \cite{GruShe87,RadWol92}; and `repetitive' in \cite{Sen95}.
In symbolic dynamical systems the term is `uniformly recurrent'.}
if for each $T > 0$ there is a finite number $M > 0$, such that
every closed ball $B\/$ of radius $M\/$ contains the center of
a translate of every possible $T$-patch in $X$. That is, for every
$T\/$-patch $\sP$ in $X$, $B\/$ contains a point of $X\/$ which is
the center of a $T$-patch of $X$ that is a translate of $\sP$.
(This $T$-patch may extend outside the ball $B$.)
The {\em repetitivity function} $M_X (T)$ is the least such radius
$M\/$ for each  $T > 0.$

Any repetitive Delone set is necessarily a Delone set of finite type,
since the definition implies that $N_X(T)$ is finite for all $T$.
The repetitivity function $M_X (T)$ of such a set provides an
additional quantitative measure of its complexity. 

A slightly different measure of repetitivity, which we denote by
$M_X'(T)$, is the least radius $M\/$ such that every ball of radius
$M\/$ contains a translate of every $T$-patch $\sP$ of $X$.
Clearly
\beql{eq110a}
 M_X'(T)=M_X(T)+T.
\eeq
Both $M_X(T)$ and $M_X'(T)$ are increasing functions of $T$.
Our reasons for choosing $M_X(T)$ instead of $M_X'(T)$ as the
fundamental repetitivity function are the following.
\begin{enumroman}
\item $M_X(T)$ is the exact analogue of the {\em recurrence function}
in symbolic dynamics, which was introduced in 1938 by Morse and Hedlund
\cite{MorHed38,MorHed40}.
\item $M_X(T)$ shares with $N_X(T)$ the property that ideal
crystals can be characterized by $M_X(T)$ being bounded for all $T$.
\item $M_X(T)$ has a natural description as the largest
of the covering radii of the sets of centers of the different
translation-equivalent $T$-patches of $X$.
\end{enumroman}

The complexity measures of a Delone set $X$ that we study in this paper
are more generally properties of the ensemble of sets that agree with $X$ 
everywhere locally.

\paragraph{Definition 1.6.}
The {\em local isomorphism class}%
\fn{The term `local indistinguishability class' has also been
proposed for this concept.}
$[[X]]$ of $X$ (under translations) consists of all sets $Y$ that
have the same atlas as $X$, i.e.
\beql{eq103}
[[X]] : = \{Y : \sA_Y = \sA_X  \}.
\eeq

Two sets, $Y$ and $Y'$, in the same local isomorphism class cannot be
distinguished by examining finite regions, because each finite region
of one appears in the other after a suitable translation. 
In measuring the complexity of sets $X$, it seems natural to  consider all
elements of a local isomorphism class $[[X]]$ to be of equal complexity. 
Many geometric quantities associated to a Delone set $X$ in $\RR^n$
are actually invariants of the local isomorphism class $[[X]]$, 
including the properties of being a finitely generated Delone set,
a Delone set of finite type, or a Meyer set. 

Local isomorphism classes of Delone sets have a natural complexity
ordering, given by the partial order $\le_{\sT}$ defined by%
\fn{The suffix $\sT$ is to indicate that the partial order
derives from inclusion of translation atlases.  In \cite{Lag98}
this partial order is considered along with a partial order
$\le_\sI$ derived from inclusion of isometry atlases.}
\beql{eq104}
[[X']] \le_\sT [[X]] \quad \mbox{if $\sA_{X'} \subseteq \sA_X$}.
\eeq
If  $f\/$ is a function whose arguments are discrete point sets
in $\RR^n$ and whose values lie in a partially ordered set,
we call $f$ {\em monotonic\/} if either
$$
[[X']]\le_\sT [[X]]\Rightarrow f(X')\le f(X)\quad\mbox{or}\quad
[[X']]\le_\sT [[X]]\Rightarrow f(X')\ge f(X),
$$
i.e.\ the ordering of $f(\cdot)$ is compatible with the partial ordering
above.  In the former case $f\/$ is {\em monotonic increasing\/} and in the
latter {\em monotonic decreasing}.  Monotonic functions are automatically 
local isomorphism class invariants. We allow such functions to take values
in the extended positive real numbers $\RR^+\cup\{0,\In\}$. 

\begin{Lemma}\label{le11} 
The following functions of a discrete point set $X$ in $\RR^n$ are
monotonic and, in particular, are local isomorphism invariants.
\begin{enumroman}
\item The Delone set parameters $r$ (monotonic decreasing) and
$R$ (monotonic decreasing).
\item The interpoint vector set $X-X$
(monotonic increasing; the partial order is inclusion).
\item The patch-counting function $N_X(T)$ (monotonic increasing).
\item The repetitivity function $M_X(T)$ (monotonic increasing).
\item The period lattice $\La_X$
(monotonic decreasing; the partial order is inclusion).
\end{enumroman}
\end{Lemma}

\paragraph{Proof.}
A function $f\/$ is monotonic increasing if it is the supremum of some
function $g\/$ defined on the atlas $\sA_X$ and is monotonic decreasing
if it is the infimum of some such function $g$.  It is straightforward
to verify that the functions in this theorem all fall into one of these
two categories.\hspace*{\fill}$\Box$\smallskip

This complexity ordering leads to a notion of `simplest' Delone sets
$X$ as those for which $[[X]]$ is a minimal element in the partial order
on local isomorphism classes.
We call such sets {\em translation-minimal Delone sets}.
In \S3 we review several equivalent notions of minimality, following
\cite{Lag98}.
For Delone sets of finite type, the notion of translation-minimality
is shown to be equivalent both to repetitivity and to minimality of
an associated topological dynamical system with an $\RR^n$-action.
The topological dynamical system  $(\sX_X, \RR^n)$ associated to
$X$ is defined analogously to tiling dynamical systems, see
\cite{RadWol92,RHB98,Rob9394,Rob,Sol97}. In the case
of a Delone set of finite type which is translation-minimal,
one has $[[X]] = \sX_X,$ while in general $[[X]] \subseteq \sX_X$.
When considering models for quasicrystals, Radin \cite{Rad91a} and
Radin and Wolff \cite{RadWol92} have argued that minimality of such
a dynamical system is analogous to being in a `ground state'.
Since the diffraction spectra of quasicrystals indicate long-range
translational order, translation-minimality seems a natural condition
to require of any structure modeling a thermodyamically stable
quasicrystalline material with zero entropy.
It is, however, a less restrictive condition, since there exist
translation-minimal sets with positive configurational entropy,
cf.\ Furstenberg~\cite[Theorem III.2]{Fur67}.

We note that the dynamical systems $( \sX_X, \RR^n)$ studied in \S3 include
systems isomorphic to any symbolic dynamical system with a $\ZZ^n$-action.
These can be encoded as Delone sets $X \subset \ZZ^n$, with translations
restricted to $\ZZ^n.$
To encode systems with a finite number of symbols, we encode the symbols as
specific patterns of points in a `block' consisting of a cubical integer
grid of a fixed size.
One can find such encodings that ensure the Delone set property and
allow the location of the edges of `blocks' to be recognized.

We consider translation-minimality to be a necessary requirement
for the `simplest' aperiodic sets.
The main object of this paper is to establish that, for a
translation-minimal Delone set of finite type, slow growth of the
repetitivity function $M_X(T)$ puts significant restrictions on
its structure, as exemplified by Theorem~\ref{th51}.
This theorem implies, in particular, that if the repetitivity
function of such a set grows at most linearly in $T$ then the
set is `diffractive' in the sense defined in \S6.

In \S4 we give the lower bound
\beql{eq108}
M_X (T) \geq r (N_X (T)^{1/n} -1 ).
\eeq
for the growth rate of the repetitivity function $M_X(T)$ in terms
of the patch-counting function $N_X(T)$ and in Corollary~4.2b
we show that there is no general upper bound, even when $X$ is
restricted to the class of Meyer sets and $N_X(T) = O( T^n)$.

The main part of the paper, \S\S5--8, concerns repetitive Delone sets that
have the slowest possible growth rate of $M_X(T)$ in one sense or another.

\paragraph{Definition 1.7.}
A repetitive Delone set of finite type is {\em densely repetitive} if
\beql{eq199}
M_X (T) = O((N_X (T))^{1/n})\quad {as} \quad T \to \In.
\eeq

By (\ref{eq108}) this is the slowest possible growth rate of $M_X(T)$
relative to the number of patches $N_X(T)$. We also formulate the
following definition bounding the growth rate of $M_X(T)$ in terms
of $T$.

\paragraph{Definition 1.8.}
A repetitive Delone set of finite type is {\em linearly repetitive} if
\beql{eq200}
M_X (T) = O(T)\quad {as} \quad T \to \In.
\eeq

Aside from ideal crystals, linearly repetitive sets have the slowest
possible growth rate of $M_X(T)$ in absolute terms, as we indicate in \S8.
However, the precise relationship between the concepts of linear and dense
repetitivity is not clear at present, and the theory would be simplified
if the following were true.

\renewcommand{\theConjecture}{\rm 1.2a}
\begin{Conjecture}
All aperiodic linearly repetitive sets are densely repetitive.
\end{Conjecture}

In \S8 we show that an equivalent form of this conjecture is the following.

\renewcommand{\theConjecture}{\rm 1.2b}
\begin{Conjecture}
Any aperiodic linearly repetitive set $X$ in $\RR^n$ satisfies
\beql{200c}
\liminf_{T \to \In}\frac{N_X(T)}{T^n} > 0.
\eeq
\end{Conjecture}
\renewcommand{\theConjecture}{\thesection.\arabic{Conjecture}}

\paragraph{Note added in proof.}  Since writing this,
Conjectures~1.2a and 1.2b have been proved by Lenz~\cite{Lenz02}.

Many constructions in the quasicrystal literature produce linearly
repetitive Delone sets.  Linear repetitivity of the Penrose tiling is
shown in Gr\"unbaum and Shephard \cite[p.563]{GruShe87}. More generally,
self-similar tiling constructions with the unique composition property
yield linearly repetitive sets, see Solomyak~\cite[Lemma 2.3]{Sol98a}.

In \S\S5--8 we prove results about the structure of densely
repetitive sets and linearly repetitive sets.
In \S5, we show that a certain class of measure-like functions 
determined by the local structure of $X\/$, which we call
`weight distributions', are approximately a constant multiple
of volume when $X\/$ is densely repetitive.
In \S6, this result is used to establish the main result of the
paper, which is that densely repetitive sets and linearly
repetitive sets have strict uniform patch frequencies.
(This notion is defined in \S6.) 
This result is equivalent to asserting that the dynamical system
with $\RR^n$-action associated to $X$ is strictly ergodic; it thus
uses a hypothesis in topological dynamics to draw a conclusion in
metrical dynamics.
The result implies that any densely or linearly repetitive set $X$ has a
uniquely defined diffraction measure in the sense of Hof~\cite{Hof95,Hof95a}.
Thus $X$ has long-range order under translations in the sense of \cite{Lag}.
We conclude \S6 with an example showing that the hypotheses of this
result are nearly best possible: there exists a repetitive Delone
set $X$ in $\RR^n$ with repetitivity function satisfying
\beql{eq111}
M_X (T) = O(T(\log T)^{2/n}(\log\log\log T)^{4/n}),
\eeq
which does not have uniform patch frequencies.

A finitely generated Delone set $X$ has an associated {\em address map}
$\phi:[X]\to\ZZ^s$, which is obtained by choosing a fixed basis of the
finitely generated Abelian group $[X]$ spanned by the vectors in $X$. 
Address maps were introduced in \cite{Lag} as a way to study the
regularity of the structure of such Delone sets.
Delone sets $X\/$ of finite type were characterized as those finitely
generated Delone sets that satisfy the Lipschitz-type condition that
there exists a constant $C$ such that
\beql{eq201}
\|\phi(\bx_1)-\phi(\bx_2)\|\le C\|\bx_1-\bx_2\|,
\quad \mbox{for }\bx_ 1,\bx_2 \in X,
\eeq
where the norms are Euclidean norms on $\RR^s$ and $\RR^n$, respectively.
Meyer sets were characterized by the existence of a linear mapping
$L: \RR^n \to \RR^s$ and a constant $C'$ such that
\beql{eq202}
\|\phi(\bx)-L(\bx)\|\le C',\quad\mbox{for }\bx\in X.
\eeq
In \S7, we show that densely repetitive sets satisfy an intermediate
property: for each densely repetitive set $X$ there exists a linear
mapping $L: \RR^n \to \RR^s$ such that
\beql{eq112}
\|\phi(\bx)-L(\bx)\|=o(\|\bx\|),
\quad \mbox{for }\bx \in X.
\eeq
For linearly repetitive sets $X$ we obtain the stronger result
that there exists a constant $\de = \de (X) > 0$ such that
\beql{JL113}
\|\phi(\bx)-L(\bx)\|=O(\|\bx\|^{1-\de})\quad\mbox{for }\bx\in X.
\eeq
We complement this by showing that there are linearly
repetitive sets $X\/$ for which there is a sequence
$\{\bx_i\}\subset X\/$ with $\| \bx_i \| \to \In$ and
$$
\|\phi(\bx_i)-L(\bx_i)\| > \|\bx_i\|^{1-\de},
$$
for some $\de<1$.

In \S8, we discuss aperiodic linearly repetitive sets as candidates
for the title of simplest aperiodic sets and propose this class
of sets as a possible notion of `perfectly ordered quasicrystal'.
We exhibit the large class of one-dimensional linearly repetitive
sets associated to Beatty sequences of the real numbers $\alpha\in[0,1]$
whose continued fraction expansions have bounded partial quotients.
Finally, we raise the question to what extent does the property of
being a linearly repetitive set $X$ put algebraic restrictions on
the image of the linear mapping $L\/$ in (\ref{eq112}).

We conclude this introduction with a brief review of related results
in the symbolic dynamics literature.
The analogue of the patch-counting function $N_X(T)$ is the
{\em block growth function}, which in higher-dimensional cases
is sometimes called the {\em (symbolic) complexity function}
$c(R)$, which counts the number of rectangular blocks of shape $R$.
It has long been known that in the one-dimensional case the bound
$c(m) \leq m$ for a single value of $m$ implies periodicity.
Coven and Hedlund~\cite{CoHe73} and Coven~\cite{Cov75} characterize
all symbolic sequences with minimal aperiodic complexity $c(m) = m + 1$.
Boshernitzan~\cite{Bos84} shows that the bound $c(m) < 3m$ implies
that the symbolic system is uniquely ergodic, hence repetitive.
Among other studies relating dynamical properties to one-dimensional
complexity are Allouche~\cite{All94} and Ferenczi~\cite{Fe96}.
A review of results concerning two-dimensional rectangular symbolic
complexity appears in Berth\'e and Vuillon~\cite{BV99}.
Various positive and negative results relating low symbolic
complexity growth to periodicity in two and higher dimensions
appear in Sander and Tijdeman~\cite{ST99a,ST99b}.
The analogue of the repetitivity function $M_X(T)$ in symbolic
dynamics might be called the {\em (uniform) recurrence function}, see
Furstenberg~\cite[p. 27]{Fur78}, and has not received much attention.

\section{Growth of the patch-counting function and periodicity}

We study the question of the extent to which slow growth of the
patch-counting function $N_X(T)$ of a Delone set of finite type
forces it to have global translation symmetries.
In \S1, we conjectured that $\limsup N_X(T)/T^n > 0$ for aperiodic
Delone sets of dimension $n$. 
Here we formulate stronger conjectures and give limited evidence for them.
In the opposite direction, for dimensions $n \geq 3$ we construct aperiodic
Delone sets of finite type that have $\liminf N_X(T)/T^n = 0$.

Sufficiently slow growth of $N_X(T)$ is known to imply the existence
of a full rank lattice of translation symmetries. 

\begin{Theorem}~\label{th21a}
If a Delone set $X$ with constants $(r, R)$ satisfies
\beql{201aa}
N_X (T)<\left\lfloor\frac{T}{R}\right\rfloor\quad\mbox{for some $T$}
\eeq
then it is an ideal crystal, i.e.\ it has a full rank lattice
$\La_X$ of translation symmetries. 
\end{Theorem}

A proof of this is given in \cite{LP}, where it is shown to be best possible,
in the sense that for every $\eps>0$ there are sets $X\/$ in $\RR^n$ with
$$
N_X(T) < \left(\frac{1}{R}+\eps\right)T\quad\mbox{for all }T>T_0
$$
that are not ideal crystals.

A version of Theorem~\ref{th21a} holds even when $N_X(T)$ is replaced
with the function $N_X^*(T)$ that counts patches up to isometry,
cf.\ Dolbilin {\em et al} \cite[Theorem 1.3]{DolLagSen}.
They showed that if $N_X^*(T) < c(n, r, R)T$ with 
\beql{207a}
c(n,r,R)^{-1}=2(n^2+1)R\log_2\Bigl(\frac{2R}{r}+2\Bigr), 
\eeq
then $X$ is an ideal crystal.
If we count patches under isometry, there is no analogue of this result
allowing a faster growth rate that forces the existence of a lower rank
lattices of periods.
Dolbilin and Pleasants \cite{DolPle} show that for each $n\ge3$ and
$\eps>0$ there is an aperiodic Delone set $X\/$ with $N_X^*(T)<T^{1+\eps}$.
The Delone sets $X\/$ in these counterexamples are not of finite type;
their construction uses twists through irrational multiples of $\pi$
to break a rank $n-1$ periodicity and make $N_X(T)$ infinite for
every $T\/$ without significantly increasing $N_X^*(T)$.
Whether such counterexamples can occur in dimension $n=2$ is an open question.

For the function $N_X(T)$, counting patches under translations, we
believe there do exist weaker conditions on the growth rate of $N_X(T)$
which will force the existence of some translation symmetries;
though none are so far known.
For definiteness, we formulate here two conjectures in arbitrary dimension,
in the strongest form that remains compatible with known counterxamples.
The first is given in the contrapositive form.

\begin{Conjecture}\label{Cj21}
For any dimension $n \geq 1$ and Delone constants $(r, R)$,
there is a positive constant $c=c(n, r, R)$ such that any Delone set
$X$ in $\RR^n$ with constants $(r, R)$ that is aperiodic satisfies
\beql{103aa}
\limsup_{T\to\In}\frac{N_X(T)}{T^n}\geq c(n,r,R).
\eeq
\end{Conjecture}

This conjecture strengthens Conjecture~\ref{Cj11} in requiring
the constant $c$ to depend only on $n$, $r$ and $R$.
If true, this result would be best possible in the sense that there
exist aperiodic Delone sets $X$ of finite type in $\RR^n$ with
$$ N_X(T) \leq (3T)^n \quad\mbox{for all }T \ge 1,$$
for example $ X := \ZZ^n\setminus\{\bf 0\}$.

Conjecture~\ref{Cj21} is the special case $j = 1$ of the following
more specific conjecture that relates slow growth of $N_X(T)$ to
the rank of the period lattice of $X$.

\begin{Conjecture}\label{Cj22}
Let the dimension $n$ be given.
For each $j$ with $1\leq j\leq n$ there is a positive constant
$c_j=c_j(n,r,R)$ such that any Delone set $X$ in $\RR^n$ with
constants $(r,R)$ that satisfies
$$\frac{N_X(T)}{T^{n+1-j}}<c_j\quad\mbox{for all }T>T_0(X)$$
has at least $j$ linearly independent periods, i.e.\ $\rank(\La_X)\geq  j$.
\end{Conjecture}

Theorem~\ref{th21a} establishes the case $j = n$ of Conjecture~\ref{Cj22}
(which includes the case $n = 1$ of Conjecture~\ref{Cj21}) in a very
strong form.  All other cases of these conjectures are open. 

We now show that Conjecture~\ref{Cj22} above cannot be strengthened to
obtain the conclusion from a single value of $T\/$, when $j\le n-2$
and $n \ge 3$.  This contrasts with the case $j = n$, where
Theorem~\ref{th21a} shows that the conclusion follows from
a single value of $T\/$ satisfying (\ref{201aa}).
In fact not even having the hypothesis for a sequence of values of
$T\/$ tending to infinity is sufficient, because we show that it is
possible for $N_X(T)/T^{n+1-j}$ to fluctuate widely as $T\to\In$.
We prove this below by modifying a construction of Sander and Tijdeman
\cite{ST99a}, who gave an example in $\RR^3$ consisting of $\ZZ^3$ with
two widely separated skew lines of lattice points removed.

\begin{Theorem}~\label{th22a}
In dimension $n\geq3$ there exist aperiodic repetitive Meyer sets $X$
such that
\beql{205a}
\liminf_{T\to\In}\frac{N_X(T)}{T^{\lceil(n+1)/2\rceil+\eps}}=0
\eeq
for every $\eps>0$.
\end{Theorem}

\paragraph{Remark.} The sets we construct satisfy
$$\limsup_{T\to\In}\frac{N_X(T)}{T^n}>0,$$
as Conjecture~\ref{Cj11} would require.

\paragraph{Proof.}
The examples we give are inductively constructed subsets of the lattice
of integer points in $\RR^n$.  We describe the case $n=3$ in detail,
but the argument generalizes easily to space of arbitrary dimension.
We start from an idea of Sander and Tijdeman \cite{ST99a} that if two skew
lines of integer points are removed from the integer lattice, the resulting
Delone set has $N(T)\le CT^2$ when $T\/$ is less than the distance between
the lines, but $N(T)>cT^3$ when $T\/$ is significantly greater than half
the distance between the lines.  The other ingredient is to notice that
if the configuration of skew lines is repeated periodically, then the
order of $N(T)$ as a function of $T\/$ declines as $T\/$ becomes large
in comparison to the fundamental region of the periodicity.
At the next stage, we alter the fundamental regions that lie along
skew lines at a larger scale.
(The simplest alteration is to replace the lines of points that were
removed from these lines of regions.)
This has the effect of again achieving $N(T)>cT^3$ when $T\/$ is
significantly greater than the distance between the new skew lines.
We now repeat this configuration periodically with a huge fundamental
region to slow down the growth of $N(T)$ again, and so on.
In this way the size of $N(T)$ can be induced to fluctuate between
$cT^3$ and not much bigger than $T^2$ as $T\/$ tends to infinity.

To make this explicit, choose a positive integer $a_1$ and remove from
$\ZZ^3$ the three systems of lines of lattice points parallel to the axes
$$
S_x=\{\bx:\bx\equiv(0,a_1,-a_1)\!\!\!\!\pmod{(1,4a_1,4a_1)}\}
$$
and the corresponding systems $S_y$ and $S_z$ with the coordinates
cyclically permuted.  (Two systems of lines would be enough---we take
three for the sake of symmetry.)  Denote by $X_1$ the resulting subset
of $\ZZ^3$ and call the lines that have been removed {\em $X_1$-lines}.
Then $X_1$ has period $4a_1$ in each of the coordinate directions.
Each $X_1$-line is the center line of an infinite right square prism of
width $2a_1$ and the interiors of these prisms are disjoint (although
the prisms of non-parallel $X_1$-lines may touch), so the minimum
distance between $X_1$-lines is $2a_1$.  Hence
$$
N_{X_1}(T)\le12T^2\quad\mbox{for}\quad T\le a_1,
$$
since any ball of radius $T\/$ meets at most one of the deleted lines
and $4T^2$ is a trivial upper bound for $N_X(T)$ when $X\/$ is $\ZZ^3$
with a line of lattice points parallel to a coordinate axis removed.

Now choose $a_2=(4b_2+1)a_1$ for some positive integer $b_2$ and let
$X_2$ be the set derived from $X_1$ by replacing the subsystem
$$
\{\bx:\bx\equiv(0,a_2,-a_2)\!\!\!\!\pmod{(1,4a_2,4a_2)}\}
$$
and the corresponding subsystems of $S_y$ and $S_z$ with the
coordinates cyclically permuted.
(Since $a_2\equiv a_1\!\!\!\!\pmod{4a_1}$, these are subsystems
of the $X_1$-lines.)
The periodicity of $X_1$ is broken on $X_2$ only by the presence
of certain lines, the {\em $X_2$-lines}, whose minimum distance
apart is $2a_2$ and consequently
$$
N_{X_2}(T)\le112a_1T^2\quad\mbox{for }a_1 < T\le a_2,
$$
since a ball of radius $T\/$ meets at most one $X_2$-line (hence
a factor $28T^2$ on the right-hand side) and the subset of $X_2$ in
the cylinder of radius $2T\/$ centered on an $X_2$-line has period
$4a_1$ in the direction of the line (hence a factor $4a_1$).
The set $X_2$ itself has period $4a_2$ in each coordinate direction.

This procedure can be iterated.  At the next step we choose
$a_3=(4b_3+1)a_2$ for some positive integer $b_3$ and form
$X_3$ by again removing points $\bx$ with
$$
\bx\equiv(0,a_3,-a_3)\!\!\!\!\pmod{(1,4a_3,4a_3)}
$$
(and their images under cyclic permutations of the axes).
For $X_3$ we have
$$
N_{X_3}(T)\le112a_2T^2\quad\mbox{for }a_2 < T\le a_3.
$$
In general we have,
\begin{equation}\label{bound}
N_{X_j}(T)\le112a_{j-1}T^2\quad\mbox{for }a_{j-1} < T\le a_j.
\end{equation}
Since the alterations needed to obtain each set from the previous
one occur further and further from the origin, the sequence of sets
$X_1,X_2,X_3,\ldots$ converges to a set $X$, which is a Meyer set
with Delone constants $r=1/2$ and $R=\sqrt5/2$.

We next note that
\begin{equation}\label{approximation}
N_X(T)=N_{X_j}(T)\quad\mbox{for }T\le a_j/2.
\end{equation}
This is because each $T$-patch $\sP$ of $X$ meets at most one $X_j$-line.
The integer points on this line are either all in $X\/$ or all in the
complement of $X$.
If their status is the same in $X\/$ as in $X_j$, then $\sP$ is also
a $T$-patch of $X_j$.
If the status in $X\/$ of the $X_j$-line is the reverse of its status
in $X_j$, then we can find a $\by\equiv\bx\!\!\!\!\pmod{4a_{j-1}}$
whose distance from all $X_j$-lines is at least $T$, and the $T$-patch
of $\by$ in $X_j$ is a translate of the $T$-patch of $\bx$ in $X$.
We also see from this argument that $X\/$ is repetitive with
$M_X(T)\le4\sqrt3a_j$ (the diameter of the fundamental region
of $X_j$) for $T\le a_j/2$.  Now
$$
N_X(T)\le112a_{j-1}T^2\quad\mbox{when }T= a_j/2,
$$
so choosing $a_j>2e^{a_{j-1}}$ gives 
$$
N_X(T)\le 112T^2\log T\quad\mbox{for }T= a_j/2,
$$
and the sequence $\{a_j/2\}$ tends to infinity.
If instead we take $T=2a_j$ and consider $T$-patches with centers in the ball
of radius $a_j-1$ about the midpoint $\bm=(-a_j,a_j,0)$ of the perpendicular
transversal of the closest $X_j$-lines to the origin in the $x\/$ and $y\/$
directions (the distance between these lines being $2a_j$), then each
$T$-patch will meet only these two $X_j$-lines in these directions,
which will be identifiable by the break in $4a_{j-1}$ periodicity, thus
allowing the position of $\bm$ within the patch to be identified.
These $T$-patches are therefore all distinct.
Hence $N_X(T)$ is of order at least $(\pi/6)T^3$ for these values
of $T$, confirming the remark preceding this proof.
In the opposite direction, it follows from (\ref{bound}) and
(\ref{approximation}) that $N_X(T)\le224T^3$ for all $T$.

The function $N_X(T)$ grows in spurts: near the $a_j$'s it increases
by a factor nearly $T\/$ while $T\/$ is increasing only by a factor~4,
but between the $a_j$'s it is almost constant.

Finally, suppose that $\bp=(p,q,r)$ is a period of $X$.
Since the planes $x=a$ that have some lines parallel to the $z$-axis deleted
all have $a\equiv a_1\!\!\!\!\pmod{4a_1}$, $p\/$ is a multiple of $4a_1$.
Since the planes $x=a$ that have a proportion less than $1/4a_1$ of the lines
parallel to the $z$-axis deleted all have $a\equiv a_2\!\!\!\!\pmod{4a_2}$,
$p\/$ is a multiple of $4a_2$.
Since the planes $x=a$ that have a proportion less than $1/4a_1$, but
greater than $1/4a_1-1/4a_2$, of the lines parallel to the $z$-axis deleted
all have $a\equiv a_3\!\!\!\!\pmod{4a_3}$, $p\/$ is a multiple of $4a_3$.
Continuing like this, we see that $p=0$.
Similarly, $q=0$ and $r=0$, so $X\/$ is aperiodic.

A similar construction works for a general dimension $n\ge3$ with
the systems of non-intersecting lines replaced by systems of
non-intersecting $\lfloor(n-1)/2\rfloor$-dimensional affine subspaces.
\hspace*{\fill}$\Box$\smallskip

\paragraph{Remark.} Theorem~\ref{th22a} implies that the hypothesis
of linear repetitivity is needed in Conjecture~1.2b, for $n \geq 3$;
repetitivity alone is not sufficient for the conclusion.

\section{Minimality}
This section formulates various notions of minimality and shows
their equivalence.
We associate to any Delone set $X$ a topological dynamical system
$( \sX_X , \RR^n )$ with an $\RR^n$-action, and show that for Delone
sets of finite type, the concepts of translation-minimality, minimality
in the sense of topological dynamics, and repetitivity are equivalent.
These results were discussed in \cite{Lag98,Sol98b}, but the topology
on the space $\sX_X$ given there is not quite right; a topology with
the required properties is given in Lenz and Stollmann~\cite{LS02}
and we follow their definition below.  We want a topology on the
collection of all closed sets in $\RR^n$ in which translation of
sets is continuous and the collection of all closed sets is compact.

We describe a topology on the collection $\sF(\RR^n)$ of closed sets
in $\RR^n$.  The {\em natural topology} on $\sF(\RR^n)$ is defined by
a neighborhood basis, as follows.  Given an integer $k \ge 1$ and
closed sets $F_1,F_2$ define the modified distance function 
$$
d_k(F_1,F_2):=\inf(\{\de:F_1\cap B({\bf 0};k)\subset F_2+B({\bf 0};\de), 
F_2\cap B({\bf 0};k)\subset F_1+B({\bf 0};\de)\}\cup \{1\}). 
$$
The function  $d_k(F_1,F_2)$ does not satisfy the triangle inequality,
but it is monotone in the cut-off parameter $k$, i.e.\ 
$d_k(F_1,F_2) \le d_{k+1}(F_1,F_2)$, and it also satisfies
$$d_k(F_1,F_2) \le d_H(F_1\cap B({\bf 0};k),F_2\cap B({\bf 0};k)),$$
where $d_H(K_1,K_2)$ is the {\em Hausdorff distance}, defined for
compact sets $K_1,K_2$ by 
$$
d_H(K_1,K_2)=\max\left(\max_{\bx_1\in K_1}\min_{\bx_2\in K_2}\|\bx_1-\bx_2\|,
\max_{\bx_2\in K_2}\min_{\bx_1\in K_1}\|\bx_1-\bx_2\|\right).
$$
The natural topology is defined by taking as a basis of neighborhoods
of a closed set $F\/$ all the sets
$$
N_{k,\eps}(F) := \{ F' \in \sF(\RR^n): d_k(F, F') < \eps\}
$$
for $k\in\ZZ^+$, $\eps>0$.  
In this topology two closed sets are `close' if they approximately
coincide in a large ball around the origin.
This topology is metrizable, but not in any canonical way, see \cite{LS02}.
It induces a topology on the set of all $(r,R)$-sets by restriction.
(The `metrics' proposed in \cite{Lag98,Sol98b} for the topology on
Delone sets are not metrics; they violate the triangle inequality.)

In particular, a sequence of Delone sets $\{X_m:m \ge 1\}$ which are all
$(r, R)$-sets converges in this topology to a limit set $X$, necessarily
also a Delone set with parameters $(r, R)$, if and only if the convergence
is `pointwise' in the sense that each point $\bx \in X$ is a limit of
a sequence of points $\bx_m \in X_m$ and, for arbitrarily large radius
$T$ and arbitrarily small $\eps>0$, for all sufficiently large $m$ all
points in $X_m$ within a distance $T$ of the origin are within $\eps$
of some point in $X$.

\paragraph{Definition 3.1.}
A {\em Delone dynamical system} $( \sX , \RR^n )$ is a collection of
Delone sets $\sX$ in $\RR^n$ that is closed in the natural topology
and also closed under translations, i.e.
\beql{eq206}
X\in\sX\Rightarrow X+\bt\in\sX \quad\mbox{for all }\bt\in\RR^n.
\eeq
It is a topological dynamical system with the $\RR^n$-action 
given by translations.

A Delone set $X$ in $\RR^n$ generates a Delone dynamical system
$(\sX_X,\RR^n)$, where $\sX_X$ is the closure in the natural topology
of the $\RR^n$-orbit of $X\/$ which is
\beql{eq207}
\sO_X : = \{X + \bt : \bt \in \RR^n \}.
\eeq
Clearly $\sO_X \subseteq [[X]] \subseteq \sX_X$.
This construction of $\sX_X$ is analogous to those developed for tiling
dynamical systems, see Radin~\cite{Rad91}, Radin and Wolff~\cite{RadWol92},
Robinson~\cite{Rob9394,Rob} and Solomyak~\cite[\S2]{Sol97}.
It gives a wide class of interesting dynamical systems, essentially
including symbolic dynamical systems as a special case%
\fn{One can encode symbols by viewing the points of the symbol
space as $k$-cubes in the lattice $k\ZZ^n$ and allowing different
configurations of points in the $k$-cubes, as mentioned in \S1.}
when the points of $X$ form a Delone subset of the lattice $\ZZ^n$. 

\paragraph{Definition 3.2.}
A Delone dynamical system $( \sX , \RR^n )$ is {\em minimal} if
$\sX_X=\sX$ for each $X\in\sX$.
Equivalently, for each $X\in\sX$ the orbit $\sO_X$ is dense in $\sX$.

\begin{Theorem}\label{th22}
Let $X$ be a Delone set of finite type in $\RR^n$.  Then:
\begin{enumroman}
\item $\sX_X$ is a compact set in the natural topology; 
\item $Y \in \sX_X \Leftrightarrow [[Y]] \le_\sT [[X]]$;
\item there exists a $Y\in\sX_X$ such that $(\sX_Y,\RR^n)$
is minimal.
\end{enumroman}
\end{Theorem}

\paragraph{Proof.}
(i) It suffices to show that every sequence of translates of
a set $X$ has a subsequence that converges to a limit Delone
set `pointwise' in the sense given above.
Convergence in this sense follows by a standard argument in
topological dynamics, using K\"onig's Infinity Lemma and the
fact that $X$ has only finitely many translation-inequivalent
patches of radius $T$, cf.\ Radin and Wolff~\cite{RadWol92},
Robinson~\cite{Rob9394,Rob}, Solomyak~\cite[Lemma~2.1]{Sol97}. 
(Alternatively, one can use the fact that $\sX_X$ is a closed subset of the
set $\sF(\RR^n)$ of all closed sets in $\RR^n$ with the natural topology,
which is compact by Lenz and Stollmann~\cite[Theorem 1.2]{LS02}.)

(ii) It is clear that if $[[Y]] \le_\sT [[X]]$ then $Y\in\sX_X$.
For the converse we have to show that every $T$-patch $\sP_Y(\by,T)$
of $Y\/$ is translation-equivalent to a $T$-patch of $X$.
Since $X\/$ is of finite type, there are only finitely many
vectors $\bv \in X-X$ with $\|\bv\| \le T$.
Let $\de$ be the minimum distance between them and put
$k = \max(\|\by\|+T+\de/2, 2/\de)$.
If $X+\bt$ is in the neighborhood $N_{k,\de/2}(Y)$ of $Y\/$ in the natural
topology then, within a radius $T+\de/2$ of $\by$, each point of either
$X+\bt$ or $Y\/$ is within a distance $\de/2$ of a unique point of the other.
In partiular, since the interpoint vector set of $Y$ is contained in that
of $X$, all interpoint vectors of corresponding points of $Y$ and $X+\bt$
within a radius $T+\de/2$ of$\by$ must match, by the definition of $\de$.
It follows that $\sP_Y(\by,T)$ is translation-equivalent to a $T$-patch of $X$.

(iii) This follows by a Zorn's lemma argument using the compactness
of $\sX_X$, see for example Furstenberg \cite[p.~29]{Fur78}.
\hspace*{\fill}$\Box$\smallskip

\begin{Theorem}\label{th23}
Let $X$ be a Delone set of finite type.
The following conditions are equivalent.
\begin{enumroman}
\item $[[X]]$ is translation-minimal. 
\item The topological dynamical system $(\sX_X,\RR^n)$ is minimal.
That is, $\sX_X = [[X]]$.
\item $X$ is repetitive. 
\end{enumroman}
\end{Theorem}

\paragraph{Proof.}
This is a form of Gottschalk's Theorem (see for example \cite[p.~136]{Pet}),
where it is proved with the group $C_\In$ in place of $\RR^n$.

(i) $\Leftrightarrow$ (ii).  This is immediate from (ii) of Theorem~\ref{th22}.

(i) $\Rightarrow$ (iii).
We prove the contrapositive.  Suppose $X$ is not repetitive.
Then there is a $T$-patch $\sP$ of $X$ for which there exists a sequence
of patches $\{\sP (\bx_i ; T_i ) : i \geq 1 \}$, with $T_i \to \In$
as $i \to \In$, none of which contain a translate of $\sP$.
Since $\sX_X$ is compact, we can extract a subsequence of the sequence
of sets $\{ X - \bx_i \}$ which converges to a limit set $Y \in \sX_X$
in the natural topology.
Since $X$ is of finite type, $Y$ contains no translate of $\sP$ and
hence $[[Y]] <_\sT [[X]]$ and $X$ is not translation-minimal.

(iii) $\Rightarrow$ (i).
Suppose that $X$ is repetitive and $[[Y]] \leq_\sT [[X]]$.
Now $Y$ contains at least one $(M_X(T)+T)$-patch, which necessarily is
a patch of $X$, so it contains a translate of every $T$-patch of $X$.
This holds for all $T$, so $\sA_X \subseteq \sA_Y$.
Hence $[[X]] = [[Y]]$ and $X$ is translation-minimal.
\hspace*{\fill}$\Box$\smallskip

\paragraph{Remark.} The compactness of $[[X]]$ is important
in modeling electronic properties of quasicrystals using
$C^*$-algebras, see Bellissard~\cite{Be92}.

\section{Repetitivity function versus patch-counting function}

In this section we study minimal Delone sets $X$ of finite type
and investigate to what extent the repetitivity function $M_X(T)$
is related to the patch-counting function $N_X(T)$.
We first derive a lower bound for $M_X(T)$ in terms of $N_X (T)$.

\begin{Theorem}\label{th31}
Let $X$ be a Delone set of finite type that is repetitive.
If $X$ is an $(r,R)$-set, then
\beql{eq401}
M_X (T) \geq r (N_X (T)^{1/n} - 1).
\eeq
\end{Theorem}

\paragraph{Proof.}
The closed ball $B({\bf 0};M_X(T))$ contains points $\bx_i$ of $X$
whose $T$-patches comprise all $N_X(T)$ translation-equivalence classes.
Since $X$ is an $(r,R)$-set, the balls of radius $r$ around these
points $\bx$ are disjoint.
They are all contained in the ball with center $\bf 0$ and radius
$M_X (T) + r$, hence
$$
\vol_n(B({\bf 0};M_X(T)+r)\geq\sum_{i=1}^{N_X (T)}\vol_n(B(\bx_i;r)).
$$
This yields
$$
\left(M_X (T)+r\right)^n\geq r^n N_X (T),
$$
from which (\ref{eq401}) follows.\hspace*{\fill}$\Box$\smallskip

It is known that there exist translation-minimal sets $X$
in one dimension that are of positive entropy.
That is, for these $X$ there is a $c>0$ such that
$$
N_X(T)\geq e^{cT}\quad\mbox{for all }T>0.
$$
Furthermore, there even exist such sets $X$ that have strict uniform
patch frequencies, as defined in \S6 (see Pleasants~\cite{Ple2}).

In the remainder of this section we show by construction
that there is no upper bound on how fast the repetitivity
function can increase, even if $X$ is required to be a 
cut-and-project set in $\RR^n$ with $N_X (T) = O(T^n )$.
Cut-and-project sets are a special subclass of Meyer sets
which have been extensively studied as models of quasicrystals,
see \cite{Lag,Moo95,Schlot,Sch,Ple1}.
The construction will be based on a result in one-dimensional
symbolic dynamics.

\paragraph{Definition 4.1.}
Let $B \in \{0,1\}^\ZZ$ be a bi-infinite 0--1 sequence.
The {\em recurrence function} $\tilde{M}_B$ of $B$ is defined by
$\tilde{M}_B (\ell)$ ($\ell = 1,2,\ldots\,$) being equal to the maximum
difference in index between the leftmost symbols of two consecutive
occurrences of any word $w \in \{0,1 \}^\ell$ of length $\ell$ that
occurs at least once in $B$.  (Here $\tilde{M}_B (\ell) = + \In$ if
some word $w$ of length $\ell$ occurs exactly once in $B$.)

\paragraph{Definition 4.2}
A {\em Beatty sequence}
$$
B_\alpha = \{ b_k \} \in \{ 0,1 \}^\ZZ
$$
is a bi-infinite 0--1 sequence associated to a real number $\alpha\in[0,1]$ by
\beql{eq405}
b_k:=\lfloor(k+1)\alpha\rfloor-\lfloor k\alpha\rfloor\quad\mbox{for }k\in\ZZ.
\eeq

All Beatty sequences for irrational $\alpha$ have $\tilde{M}_{B_\alpha}(k)$
finite for all $k \geq 1$, see \cite{MorHed40}.
We call a function $g : \RR^+ \to \RR^+$  a {\em growth function}
if it is continuous and nondecreasing.

\begin{Theorem}\label{th32}
Given any growth function $g$ there exists an irrational $\alpha\in[0,1]$
whose Beatty sequence $B = B_\alpha$ has recurrence function
\beql{eq406}
\tilde{M}_B(\ell)>g(\ell)\quad\mbox{for infinitely many $\ell\geq 1$.}
\eeq
\end{Theorem}

\paragraph{Proof.}
Morse and Hedlund \cite[Theorem 9.1 and Lemma 10.3]{MorHed40},
show that the recurrence function of $B$ is explicitly given by
\beql{eq407}
\tilde{M}_B(\ell)=q_k+q_{k+1}\quad\mbox{for }q_k\leq\ell<q_{k+1},
\eeq
where $p_k/q_k$ is the $k\/$th convergent of the regular continued
fraction expansion of $\alpha$.  If we take $\alpha$ to be the number
whose continued fraction expansion is $[0,a_1,a_2,a_3,\ldots ]$,
where the partial quotients $a_1,a_2,\ldots$ are chosen successively
to satisfy $a_{k+1} \ge g(q_k)$ for $k \ge 0$, then we have
$$
\tilde{M}_B(q_k)=q_k+q_{k+1}=q_k+a_{k+1}q_k+q_{k-1}>a_{k+1}\ge g(q_k)
$$
for $k \ge 0$.\hspace*{\fill}$\Box$\smallskip

\renewcommand{\theCorollary}{\rm 4.2a.}
\begin{Corollary}
Given any growth function $g(T)$ there exists a translation-minimal 
one-dimensional cut-and-project set $X$ projected from $\ZZ^2$
using $[0,1)$ as a window whose repetitivity function satisfies
\beql{eq408}
M_X (T_i) > g (T_i),i = 1,2,\ldots\,,
\eeq
for a sequence of $T_i \to \In$.
\end{Corollary}

\paragraph{Proof.}
For any real number $\alpha \in (0,1)$, the one-dimensional
cut-and-project set $X\/$ obtained by projecting the integer
lattice in $\RR^2$ orthogonally on to a line of slope $\alpha$,
using the interval $[0,1)$ of the $y$-axis as the window, is
identical to the sequence of intervals obtained by replacing
the symbols 0 and 1 in the Beatty sequence $B=B_\alpha$ of
$\alpha$ by intervals of lengths $s\/$ and $\ell\/$, where
$s\/$ and $\ell\/$ depend on $\alpha$ and satisfy
$1/\sqrt2 < s < \ell < \sqrt2$.  Consequently
$$
\hspace{54.5mm}M_X(T)>\tilde{M}_{B_\alpha}(\lfloor T/\sqrt2\rfloor).
\hspace{51.5mm}\Box
$$
\smallskip

\renewcommand{\theCorollary}{\rm 4.2b.}
\begin{Corollary}
Given any growth function $g(T)$, there is a translation-minimal
aperiodic Meyer set $Y$ in $\RR^n$ with patch-counting function
$N_Y(T)=O(T^n)$ whose repetitivity function satisfies
\beql{eq409}
M_Y (T_i) > g(T_i)\quad i = 1,2,\ldots\,,
\eeq
for a sequence of $T_i \to \In$.
\end{Corollary}
\renewcommand{\theCorollary}{\thesection.\arabic{Corollary}}

\paragraph{Proof.}
All Beatty sequences $B = B_\alpha$ with irrational $\alpha$ are
aperiodic and have
$$
N_B (k)=k+1\quad\mbox{for all }k=1,2,\ldots.
$$
Thus the cut-and-project sets $X$ constructed in Corollary~4.2a
are aperiodic and have $N_X (T) = O(T)$.
Take $Y=X^n\subseteq\RR^n$ to be the direct product of $n$ such sets.
Since every cut-and-project set $X\/$ is a Meyer set 
(\cite[Corollary 5.6]{Moo97}) and since direct products
of Meyer sets are Meyer sets, $Y\/$ is a Meyer set.
Clearly $N_Y (T) = O (T^n )$ and $M_Y(T) \geq M_X(T)$,
so (\ref{eq409}) is satisfied.\hspace*{\fill}$\Box$\smallskip

The Meyer sets constructed in Corollary~4.2b are cut-and-project sets
from $\RR^{2n}$ to $\RR^n$ with the window being the hypercube $[0,1)^n$.

\section{Locally defined functions on repetitive sets}

In this section we show that the properties of linear repetitivity and
dense repetitivity of a Delone set $X\/$ in $\RR^n$ cause a large class
of functions defined locally with reference to $X\/$ to have uniquely
determined average values on $\RR^n$.
These general results will be applied in \S6 to patch frequencies
(with implications for autocorrelation measures and diffraction) and in
\S7 to estimating path lengths in address space in the course of showing
that, for both densely repetitive and linearly repetitive sets, the
address map can be reasonably well approximated by a linear function.

We shall consider functions defined on `boxes' (that is, direct
products of coordinate intervals) in $\RR^n$ and express our results
in terms of the volume, $\vol_n (B)$, surface area, $\sig(B)$, and
width, $\om(B)$, of a box $B$.

\paragraph{Definition~5.1.}
A {\em box} in $\RR^n$ is a subset of the form
$$
B=\{(x_1,x_2,\ldots,x_n)\mid a_i\le x_i\le b_i\ (i=1,\ldots,n)\},
$$
where $a_i<b_i\in\RR$ for each $i$.  If $l_i=b_i-a_i$ ($i=1\ldots n$), then
the {\em volume, surface area\/} and {\em width\/} of $B\/$ are given by
$$
\vol_n (B)=\prod_{i=1}^nl_i,
\qquad\sig(B)=2\,\vol_n (B)\sum_{i=1}^n{1\over l_i},
\qquad\om(B)=\min_{1\le i\le n}l_i.
$$

The class of functions we deal with is the following.

\paragraph{Definition~5.2.}
A {\em local weight distribution} on a Delone set $X\/$ in $\RR^n$
is a function $\bw$, with values in $\RR^s$ whose domain is the set
of all boxes $B\/$ in $\RR^n$ with $\om(B)>U_0$ (for some $U_0\ge0$)
and which has the following properties:
\begin{enumalph}
\item $\bw(B)=O(\vol_n (B))$;
\item $\bw$ is `approximately translation-invariant' on $X$,
in the sense that if
$$
(B+\bt)\cap X=(B\cap X)+\bt
$$
for some $\bt\in\RR^n$, then
$$
\bw(B+\bt)=\bw(B) + O(\sig(B));
$$
and
\item $\bw$ is `approximately additive', in the sense that
if $B=B_1\cup B_2\cup\ldots\cup B_k$ is a partition of $B\/$
into non-overlapping smaller boxes, then
$$
\bw(B)=\sum_{i=1}^k \bw(B_i)+O\left(\sum_{i=1}^k\sig(B_i)\right).
$$
\end{enumalph}
In this definition the $O$-constants depend only on $\bw$.

Let $\sB(U)$ denote the set of `squarish' boxes, all of whose side
lengths $l_i$ lie in the range $U\le l_i\le 2U$.
Given a real-valued local weight distribution $w\/$ on a Delone set
$X$, we define, for $U>U_0$, the upper and lower local densities
\beql{eq501}
f^+(U)=\sup_{B\in\sB(U)}w(B)/\vol_n (B)\quad\mbox{and}\quad
f^-(U)=\inf_{B\in\sB(U)}w(B)/\vol_n (B).
\eeq

\begin{Lemma}\label{Le41}
If $w\/$ is a real-valued weight distribution on a Delone set $X\/$ 
and the upper and lower densities $f^+$, $f^-$ are as in (\ref{eq501})
then $f^+(U)$ and $f^-(U)$ tend to limits as $U\to\In$.
\end{Lemma}

\paragraph{Proof.} 
We show that $f^+$ is an `approximately decreasing' function
of $U\/$ (and $f^-$ is `approximately increasing').

If $W\ge U$, then any box $B\in\sB(W)$ can be subdivided into boxes
$B_i$ of $\sB(U)$.  (This is clearly possible in the case $n=1$ of an
interval on the line and in higher-dimensional cases a subdivision
that is a direct product of one-dimensional subdivisions can be used.)
Since $\sum \vol_n (B_i)=\vol_n (B)$, the number of boxes $B_i\in\sB(U)$
needed to cover $B\/$ is $O(W^n/U^n)$.  Also, $\sig(B_i)=O(U^{n-1})$ for
each of them.  So Definition~5.2(c) (on division by $\vol_n (B)$) gives
\beql{eq502}
f^+(W)<f^+(U)+{C_1\over U},
\eeq
for some constant $C_1$.

In view of (a) of Definition~5.2, $\liminf f^+$ is finite.
For all $\eps>0$, there is a $U>1/\eps$ with
$$
f^+(U)<\liminf_{U\to\In}f^+(U)+\eps,
$$
so by (\ref{eq502})
$$
f^+(W)<\liminf_{U\to\In}f^+(U)+(1+C_1)\eps\quad\mbox{for all }W>U,
$$
and letting $W\to\In$ gives
$$
\limsup_{W\to\In}f^+(W)\le\liminf_{U\to\In}f^+(U)+(1+C_1)\eps.
$$
Since $\eps$ can be chosen arbitrarily small we have
$$
\limsup f^+\le \liminf f^+,
$$
showing that $f^+(U)$ tends to a limit as $U\to\In$.

The proof that $f^-(U)$ tends to a limit is similar.
\hspace*{\fill}$\Box$\smallskip

\begin{Theorem}\label{th41}
If $X\/$ is a linearly repetitive Delone set then there is a
constant $\de=\de(X)>0$ such that every local weight distribution
$\bw$ on $X\/$ has an average value $\bff$ in $\RR^s$ satisfying
\beql{JL45}
\frac{\bw(B)}{\vol_n (B)}=\bff+O(\om(B)^{-\de})
\eeq
as $\om(B)$ tends to infinity.
(The $O$-constant may depend on $\bw$ as well as $X$.)
\end{Theorem}

\paragraph{Proof.} 
There is no loss of generality in assuming that $\bw = w\/$
is real-valued, since a vector-valued $\bw$ can be treated by
applying the result for real-valued $w\/$ to each coordinate.
For real-valued $w$, the existence of the constant field $\bff = f\/$
is equivalent to the two limits of Lemma~\ref{Le41} being the same.
We denote these limits by $f^+$ and $f^-$.

Throughout the proof we shall always take $U\ge U_1=\max(U_0,R,1)$.
For such a $U\/$ we can find a box $B_1'\in\sB(U)$ with
$$
\frac{w(B_1')}{\vol_n (B_1')}<f^-(U)+\frac{1}{U}.
$$
Now $B_1'$ is contained in some $V$-patch of $X$, where
$V=(2\sqrt n+1)U>2\sqrt n\,U+R$, and by linear repetitivity every ball of
radius $C'V\/$ contains a point of $X\/$ whose $V$-patch is a translate of
this one, where $C'$ is the implied constant in (\ref{eq200}).  Consequently,
if $W=2(C'V+V+U)=C_2U$, where $C_2=2(C'+1)(2\sqrt n+1)+2$, then every box
$B\in\sB(W)$ contains a box $B_1$ which is translation-equivalent to
$B_1'$ on $X\/$ and whose distance from the boundary of $B\/$ is $\ge U$.
In a similar way as in the proof of Lemma~\ref{Le41}, we can subdivide
$B\/$ into boxes of $\sB(U)$, one of which is $B_1$.  Also
$$
\vol_n (B_1)\ge U^n=C_2^{-n}W^n\ge c\,\vol_n (B)
$$
with $c=(2C_2)^{-n}>0$.
So, by the reasoning that led to (\ref{eq502}), we obtain
\beql{eq503}
f^+(W)<c\left(f^-(U)+\frac{1}{U}\right)+(1-c)f^+(U)+\frac{C_1}{U}
\eeq
and letting $U\to\In$ gives
$$
c f^+\le c f^-,
$$
showing that $f^+=f^-=f$, say.

To establish the convergence rate, we first repeat the above argument
to obtain
\beql{eq504}
f^-(W) > c f^+(U) + (1-c)f^-(U) - \frac{C_1+1}{U}
\eeq
then subtract (\ref{eq504}) from (\ref{eq503}) to obtain
\beql{eq505}
\De(W) < (1-2c)\De(U) + \frac{2(C_1+1)}{U},
\eeq
where $\De(U)=f^+(U)-f^-(U)\ge0$.

We now show that
$$
\De(U) < C_3 U^{-\de}\quad\mbox{for }U\ge U_1,
$$
where
\beql{eq506}
\de=\log(1-c)^{-1}/\log C_2
\eeq
and $C_3$ is chosen large enough to ensure that $\De(U) < C_3 U^{-\de}$,
for $U_1\le U\le C_2U_1$, and $c C_3 > 2(C_1 +1)$.
Assume, as an induction hypothesis, that
$$
\De(U) < C_3 U^{-\de}\quad\mbox{for }U_1\le U\le C_2^kU_1.
$$
This holds for $k=1$.  Given $W\le C_2^{k+1}U_1$, we have
\begin{eqnarray*}
\De(W)
&<&(1-2c)\De(W/C_2)+\frac{2(C_1+1)C_2}{W},\quad\mbox{by (\ref{eq505}),}\\
&<&C_3(1-c)W^{-\de}C_2^\de,\quad\mbox{by the induction hypothesis,}\\
&=&C_3 W^{-\de},\quad\mbox{by (\ref{eq506}).}
\end{eqnarray*}
This completes the induction and shows that $\De(U)<C_3 U^{-\de}$
for all $U>U_1$.

Now letting $W\to\In$ in (\ref{eq502}) (and its analog for $f^-$) gives
$$
f^-(U)-\frac{C_1}{U}<f<f^+(U)+\frac{C_1}{U}
$$
so
\beql{eq507}
\left|\frac{w(B)}{\vol_n(B)}-f\right|\le\De(U)+\frac{C_1}{U}=O(U^{-\de})
\eeq
for $B\in\sB(U)$, since $\de<1$.

Finally, an arbitrary box $B\/$ with $\om(B)>U_1$ can be subdivided into
$O(\vol_n (B)/\om(B)^n)$ boxes of $\sB(\om(B))$ with total surface area
$O(\vol_n (B)/\om(B))$, so (\ref{eq507}) and (c) of Definition~5.2 give
$$
\frac{w(B)}{\vol_n(B)}=f+O(\om(B)^{-\de}),
$$
completing the proof.\hspace*{\fill}$\Box$\smallskip

\begin{Theorem}\label{th42}
If $X\/$ is a densely repetitive Delone set and $\bw$ is a local
weight distribution on $X\/$, then there is an $\bff$ in $\RR^s$ with
$$
\frac{\bw(B)}{\vol_n(B)}=\bff+o(\om(B)),
$$
as $\om(B)$ tends to infinity.
\end{Theorem}

\paragraph{Proof.} 
Again there is no loss of generality in assuming that $\bw=w$
is real-valued.

Given $U>U_0$ we have defined $f^+(U)$ and $f^-(U)$.
We now also define $f^0(U)$ as follows.
Let $\sP_X(\bx_1,\sqrt n\,U),\ldots,\sP_X(\bx_N,\sqrt n\,U)$
be a complete set of $N=N(\sqrt n\,U)$ translation-inequivalent
$\sqrt n\,U$-patches of $X\/$ and let
$$
w_i = w_i(U) = w(\sC(\bx_i,U))\quad(i=1,\ldots,N),
$$
where $\sC(\bx_i,U)$ is the cube of side $2U\/$ with center
$\bx_i$ (which is contained in $\sP_X(\bx_i,\sqrt n\,U)$).
We take $f^0(U)$ to be the median of the set of numbers
$$
\{ w_i/(2U)^n : i = 1,\ldots,N(\sqrt n\,U) \},
$$
that is, the middle number when the set is in size order (when $N\/$ is odd)
or the average of the two middle numbers (when $N\/$ is even).  There is an
ambiguity in this definition as it stands, because two points $\bx_i^{}$,
$\bx_i'$ with translation-equivalent $\sqrt n\,U$-patches may nevertheless
give different values of $w_i$.  However, by Definition~5.2(b) the values
differ by at most a constant multiple of $U^{n-1}$ and, consequently,
$f^0(U)$ is defined to within a constant multiple of $1/U\/$ and its
behavior as $U\/$ tends to infinity is unaffected by the ambiguity.
The definition could be made unambiguous by adopting some specific
rule for choosing the points $\bx_1,\ldots,\bx_N$.

\begin{sloppypar}
We now choose a sequence $U_1,U_2,\ldots\,$, tending to infinity, with
$U_1 > U_0$ and $M_X(\sqrt n\,U_k)\ge \sqrt n\,U_k$ for $k = 1,2,\ldots$.
(If there is no such infinite sequence then $X\/$ is linearly repetitive
and the result is a consequence of Theorem~\ref{th41}.)
Put $W_k=7M_X(\sqrt n\,U_k)$ and consider an arbitrary box $B\in\sB(W_k)$.
To get a non-trivial bound for $w(B)$ it is no longer enough to find a
single sub-box $B_1$ of $B\/$ in $\sB(U_k)$ with $w(B_1)$ small since,
when the growth of $M_X$ is faster than linear, this no longer fills
a positive proportion of $B$.  We need a large number of well-spaced
sub-boxes with small $w$.  By the definition of $f^0(U_k)$, there are
$\ge N/2 = N(\sqrt n\,U_k)/2$ translation-inequivalent $\sqrt n\,U_k$-patches
of $X\/$ such that if $\sP$ is a patch translation-equivalent to any one
of them and $\sC$ is the coordinate cube of side $2U_k$ inscribed in $\sP$,
then $w(\sC)/(2U_k)^n \le f^0(U_k)+C_4/U_k$ for some $C_4$.  By the
definition of $M_X$, every cube of side $2M_X(\sqrt n\,U_k)+\sqrt n\,U_k$
contains at least one patch translation-equivalent to each of these.
Hence $B\/$ contains at least $N/2$ subcubes $\sC(\bx,U_k)$ of side $2U_k$
for each of which $w/(2U_k)^n\le f^0(U_k)+C_4/U_k$ and none of which comes
within a distance $2\sqrt n\,U_k$ of the boundary of $B$.  We require a
set of subcubes that do not overlap and, in fact, it will be convenient to
have them a distance at least $4\sqrt n\,U_k$ apart.  Since each of our
subcubes $\sC(\bx,U_k)$ has a point $\bx$ of $X\/$ at its center and a
ball of radius $6\sqrt n\,U_k$ contains at most $((6\sqrt n\,U_k+r)/r)^n$
points of $X$, each subcube $\sC(\bx,U_k)$ we pick excludes at most
$((6\sqrt n\,U_k+r)/r)^n-1$ others.  Consequently, we can choose at least
$\frac{1}{2}r^n N/(7\sqrt n\,U_k)^n\ge\frac{1}{2}(r M/7\sqrt n\,C''U_k)^n$
of these subcubes of $B\/$ that are at least a distance $4\sqrt n\,U_k$
apart, where $C''$ is the implied constant in (\ref{eq199}) and
$M = M_X(\sqrt n\,U_k)$.
\end{sloppypar}

We next verify that $B$ has a partition into boxes of $\sB(U_k)$ which
includes all these subcubes.  This can be done by thinking of $\RR^n$ as
composed of a fixed background grid of packed cubes of side $U_k$.  Each
chosen cube $\sC(\bx,U_k)$ meets a set of grid cubes forming a box-shaped
block, which we enlarge by the inclusion of an outer single-layer shell
of grid cubes disjoint from $\sC(\bx,U_k)$.  Because of the separation
of the chosen cubes, their enlarged blocks of grid cubes do not overlap.
As in previous proofs, this enlarged block can be partitioned into boxes
of $\sB(U_k)$ one of which is $\sC(\bx,U_k)$.  Similarly, we consider the
two-layer shell of grid cubes consisting of those that meet the boundary
of $B$ together with an inner single-layer lining shell.  Again, the part of
this two-layer shell within $B$ can be partitioned into boxes of $\sB(U_k)$
and does not overlap the shells of the chosen cubes because of the distance
of these cubes from the boundary of $B$.  So $B$ can be partitioned into
the cubes $\sC(\bx,U_k)$ (of side $2U_k$), the boxes of $\sB(U_k)$ just
described and the remaining grid cubes (of side $U_k$).

\sloppy{Since $\vol_n(B)\le(14M)^n$ and the total volume of the cubes
$\sC(\bx,U_k)$ is at least $\frac{1}{2}(2rM/7\sqrt n\,C'')^n$,
Definition~5.2(c), on division by $\vol_n (B)$, gives}
$$
f^+(W_k)<c\left(f^0(U_k)+\frac{C_4}{U_k}\right)+(1-c)f^+(U_k)+\frac{C_1}{U_k},
$$
where $c=\frac{1}{2}(r/49\sqrt n\,C'')^n$.  On letting $k\to\In$ we obtain
$$
\liminf_{U\to\In}f^0(U)\ge f^+.
$$
A similar argument, choosing cubes $\sC \subseteq B$ with
$w/(2U_k)^n \ge f^0(U_k) - C_4/U_k$, leads to
$$
\limsup_{U\to\In}f^0(U)\le f^-.
$$
Hence $f^+ = f^-$.  The proof is completed as in the last two paragraphs
of the proof of Theorem~\ref{th41}.\hspace*{\fill}$\Box$\smallskip

We end this section by showing that linear repetitivity is not the exact
cut-off point for the universal existence of unique averages for local
weight distributions.  It is possible to trade some of the error term
in Theorem~\ref{th41} for a weakening of the hypothesis and show that
every marginally superlinearly repetitive Delone set has this property.

\begin{Theorem}\label{th43}
If $X\/$ is a Delone set in $\RR^n$ satisfying
$$
M_X(T) = O(T(\log T)^{1/n})
$$
as $T\/$ tends to infinity and $\bw$ is a local weight
distribution on $X\/$ then there is an $\bff$ in $\RR^s$ with
$$
\frac{\bw(B)}{\vol_n (B)}=\bff+o(\om(B)),
$$
as $\om(B)$ tends to infinity.
\end{Theorem}

\paragraph{Proof.}
We first note that the derivation of (\ref{eq503}), (\ref{eq504}) and
(\ref{eq505}) in the proof of Theorem~\ref{th41} remains valid when $C'$
is a function of $T\/$ instead of being constant.  If we assume that 
$$
M_X(T) < CT(\log T)^{1/n} \quad \mbox{for $T \ge U'$,}
$$
where we can suppose that $U' \ge \max (U_0,2\sqrt n+1)$, and take
$W = C_5 U (\log U)^{1/n}$ with $C_5 = 4(2\sqrt n +1)C +2$, then
(\ref{eq505}) becomes
$$
\De(W)<\left(1-\frac{2c_1}{\log U}\right)\De(U)+\frac{2(C_1+1)}{U},
$$
where $c_1 = (2C_5)^{-n}$ and $C_1$ is as before.
If we put $E(T)=\max(\De(T),T^{-1/2})$ then we can find a $U_1$ such that
\beql{eq508}
E(W)<\left(1-\frac{c_1}{\log U}\right)E(U)\quad\mbox{for $U\ge U_1$.}
\eeq
We can also take $U_1 \ge \max(C_5,4)$.

We now define an increasing sequence $U_1 < U_2 < \cdots$ by
$$
U_{k+1} = C_5 U_k (\log U_k)^{1/n}.
$$
Clearly $U_k \to \In$ as $k \to \In$.
To bound $U_k$ above, we use an induction argument to show that
$$
\log U_k < U_1 k \log (U_1 k) \log \log (U_1 k).
$$
For the induction step we have
\begin{eqnarray*}
\log U_{k+1}&=&\log U_k+\frac{1}{n}\log\log U_k+\log C_5\\
&<&U_1k\log(U_1k)\log\log(U_1k)+\log(U_1k)+\log\log(U_1k)\\
&&+\log\log\log(U_1k)+\log C_5,\quad\mbox{by the induction hypothesis,}\\
&<&U_1(k+1)\log(U_1k)\log\log(U_1k),
\quad\mbox{because $U_1\ge4$ and $U_1\ge C_5$,}\\
&<&U_1(k+1)\log(U_1(k+1))\log\log(U_1(k+1)).
\end{eqnarray*}

Now, for any $U \ge U_1$, let $K\/$ be the largest integer with
$U_K \le U$.  Then repeated application of (\ref{eq508}) gives
$$
E(U)<C_6\prod_{k=2}^K
\left(1-\frac{c_1}{U_1k\log(U_1k)\log\log(U_1k)}\right),
$$
where $C_6$ is an upper bound for $E(U)$ in the range $U_1\le U<U_2$.
Since the product on the right-hand side diverges to zero as $K\/$
tends to infinity, $E(U)$ (and hence $\De(U)$) tends to zero as $U\/$
tends to infinity.  The proof is now completed as in the last two
paragraphs of the proof of Theorem~\ref{th41}.\hspace*{\fill}$\Box$

\section{Repetitive sets and uniform patch frequencies}

One form of long range order in a discrete set $X$ is the property
that, for each $T$, all of its $T$-patches are evenly distributed
in the sense of having limiting frequencies.  We show that densely
repetitive sets and linearly repetitive sets have this property.

\paragraph{Definition 6.1.}
Given a $T$-patch $\sP$ of $X$ translated to {\bf 0},
and a bounded region $D\subset\RR^n$ define
$$
n_{\sP}(D) := |\{ \bx\in D : \sP+\bx \subseteq X \}|.
$$
That is, $n_{\sP}(D)$ counts the number of points of $X\/$ in
$D\/$ whose $T$-patches are translation equivalent to $\sP$.  

\paragraph{Definition 6.2.}
A Delone set $X$ has {\em uniform patch frequencies} if every patch
$\sP$ has a uniform limiting frequency $f( \sP )$ in the sense that
\beql{eq601}
\frac{n_\sP(B(\bt;U))}{\kappa_n U^n} = f(\sP) + o(1)
\eeq
uniformly in $\bt$ as $U\to\In$, where $\kappa_n$ is the volume of
the unit ball in $\RR^n$ (so that $\kappa_n U^n=\vol_n(B(\bt;U))$).
It has {\em strict uniform patch frequencies} if, in addition,
$f(\sP)$ is non-zero for every patch $\sP$ that occurs in $X$.

A topological dynamical system on a compact space $\Omega$ with an
$\RR^n$-action is {\em uniquely ergodic} if it has a unique probability
measure on Borel sets that is invariant under the $\RR^n$-action.
The property of having uniform patch frequencies is equivalent to
the dynamical system $(\sX_X,\RR^n)$ being uniquely ergodic, see
Lee, Moody and Solomyak~\cite[Theorem 2.7]{LMS02}.
The extra property that the frequency of every patch that occurs in $X$ is
positive implies that $(\sX_X,\RR^n)$ is minimal and, hence, strictly ergodic.
Lunnon and Pleasants~\cite{LP1} take the property of having strict
uniform patch frequencies as part of the definition of being strongly
quasicrystallographic.

The main result of this paper is as follows. 
\begin{Theorem}\label{th51}
If $X$ is a Delone set in $\RR^n$ that is either linearly repetitive
or densely repetitive, then $X$ has strict uniform patch frequencies.
In the case that $X$ is linearly repetitive the error term in (\ref{eq601})
has the form $O(U^{-\de})$, uniformly in $\bt$, where $\de>0$ depends
only on $X$ but the $O\/$-constant may depend on the patch $\sP$. 
\end{Theorem}

\paragraph{Proof}
We note that $n_{\sP}$ is a local weight distribution on $X$.  It clearly
satisfies (a) of Definition~5.2 and (c) (with no error term necessary).
The error term in (b) arises from points $\bx$ in $B+\bt$ whose
$T$-patches extend outside $B+\bt$.  Such points are confined to
a border region of $B+\bt$ with thickness $T$, whose volume is at
most $T\sig(B)$, and the number of them is $O(\sig(B))$.

Now Theorems~\ref{th41} and \ref{th42} ensure the existence of
a limiting frequency $f(\sP)$ such that
\beql{eq603}
\frac{n_{\sP}(B)}{\vol_n(B)} = f(\sP) + \eps(B)
\eeq
with $\eps(B)\to 0$ as $\om(B)\to\In$ in both the linearly and
densely repetitive cases.  In the linearly repetitive case
$\eps(B) = O(\om(B)^{-\de_1})$, where $\de_1$ is the $\de$
of Theorem~\ref{th41}.
To see that $f(\sP)>0$, note that for every large integer $N\/$,
every cube of side $2M_X(T)N\/$ can be subdivided into $N^n$ cubes
of side $2M_X(T)$, so contains at least $N^n$ distinct points of
$X\/$ whose $T$-patches are translation equivalent to $\sP$.
Letting $N\to\In$ now gives $f(\sP) \ge (2M_X(T))^{-n} > 0$.

To complete the proof we need to make the transition from boxes to balls.
Given a ball $B(\bt; U)$, we can cover it with $k=O(U^{n/2})$ cubes
$C_1,\ldots,C_k$ of side $\sqrt U$ of which at most $O(U^{(n-1)/2})$
are not contained in $B(\bt;U)$.  So
$$
\vol_n(B(\bt;U)) = kU^{n/2} - O(U^{n-1/2})
$$
and the error in estimating $n_{\sP}(B(\bt;U))$ as $\sum n_{\sP}(C_i)$
is also $O(U^{n-1/2})$.
These approximations to the volume and patch-count of $B(\bt;U)$ lead to
\begin{eqnarray*}
\frac{n_{\sP}(B(\bt;U))}{\vol_n(B(\bt;U))}
&=&\frac{1}{k U^{n/2}}\sum_{i=1}^k n_{\sP}(C_i) + O(U^{-1/2})\\
&=&f(\sP) + \frac{1}{k}\sum_{i=1}^k \eps(C_i) + O(U^{-1/2})
\end{eqnarray*}
by (\ref{eq603}).  The last two terms on the right-hand side tend to zero
as $U\to\In$ in both the linearly and densely repetitive cases and in the
linearly repetitive case they are $O(U^{-\de})$, with $\de=\de_1/2$.
\hspace*{\fill}$\Box$\smallskip

The major consequence of Theorem~\ref{th51} is that $X\/$ is
diffractive in the sense of Hof~\cite{Hof95}, as we now explain.

\paragraph{Definition 6.3.}
An {\em autocorrelation measure} (in the sense of Hof~\cite{Hof95,Hof95a})
of a discrete set $X$ in $\RR^n$ is any limit measure in the vague topology,
as $T\to\In$, of the set of discrete measures $\mu_T$ defined for $T>0$ by
\beql{eq602}
\mu_T:=\frac{1}{\kappa_n T^n}\sum_{\|\bx_1\|,\|\bx_2\|<T}\de_{\bx_1-\bx_2}.
\eeq

In this definition, we view a measure $\mu$ as a linear functional on the
space $\sK$ of complex-valued continuous functions $f\/$ with compact
support on $\RR^n$ which has the additional property that, for any
compact set $E$, there is a constant $a_E$ (depending on $\mu$) such that
$$
|\mu (f)|\leq a_E \|f\|\quad\mbox{if $\supp (f)\subseteq E$},
$$
where $\|f\|=\sup\{|f(\bx)|:\bx\in\RR^n\}$, and we say that a sequence of
measures $\{\mu_k\}$ {\em converges to $\mu$ in the vague topology} if
$$
\lim_{k\to\In}\mu_k(f)=\mu(f)\quad\mbox{for all $f\in\sK$}.
$$

Any Delone set $X$ has at least one autocorrelation measure%
\fn{This follows by the remark before \cite[Prop.~2.2]{Hof95}
because the measure $\mu_X : = \sum_{\bx \in X} \de_\bx$
is translation-bounded, as defined in \cite{Hof95}.}
and typically (when $X\/$ has large scale irregularities) has
many autocorrelation measures.

\paragraph{Definition 6.4.}
A Delone set $X$ is {\em diffractive} if it has a unique autocorrelation
measure $\gamma_X$. The {\em diffraction measure} of $X$ is the Fourier
transform $\hat{\gamma}_X$ of the autocorrelation measure $\gamma_X$. 

The diffraction measure $\hat{\gamma}$ is a mathematical analog of the
far-field x-ray diffraction spectrum of $X$, see Hof~\cite{Hof95,Hof95a}.
It is necessarily a positive measure. Note that this definition of
a Delone set being diffractive does not require that its diffraction
measure contain any pure point component.

The existence of uniform patch frequencies implies that $X$ has a
unique autocorrelation measure $\gamma = \gamma_X$ and hence that
it is diffractive in the above sense.  It also implies that $X$
has unique $k$-point correlation measures for all $ k \geq 1$.

Theorem~\ref{th51} now yields the following.

\begin{Corollary}\label{Co51}
If $X$ is a linearly or densely repetitive Delone set in $\RR^n$,
then it has a unique autocorrelation measure $\gamma_X$.
This measure $\gamma_X$ is a pure discrete measure supported on $X-X$.
In particular $X$ is diffractive.
\end{Corollary}

\paragraph{Proof.}
The existence and uniqueness of $\gamma_X$ follows from Example~2.1 in
Hof~\cite{Hof95}.  Now $X-X$ is a uniformly discrete set, since $X$ is
a Delone set of finite type, hence all measures $\mu_T$ in (\ref{eq602})
are supported on $X-X$ and $\gamma_X$ inherits this property.
\hspace*{\fill}$\Box$\smallskip

Theorem~\ref{th51} also applies to symbolic dynamical systems given by a
two-sided shift $\Sigma$ on a finite alphabet $\sA = \{0, 1,...,m - 1\}$.
Elements of such a shift can be encoded as Delone sets of finite type on
the line $\RR$, by assigning $m\/$ intervals of different fixed lengths
to the symbols, and the shift is minimal with a recurrence function that
grows at most linearly in the word length $\ell\/$ if and only if the
associated Delone sets are linearly repetitive.
Clearly such linearly repetitive shifts have the number of words
of length $\ell\/$ bounded above by $C\ell\/$ for some constant
$C\/$ (by the analog for shifts of Theorem~\ref{th31}).
Theorem~\ref{th51} implies that linearly repetitive shifts are
uniquely ergodic.  Boshernitzan~\cite{Bos84} proves a related result:
any minimal shift whose word-count function is bounded above by
$C\ell\/$ has finitely many ergodic invariant measures, the number
of them bounded in terms of $C\/$ and the number of symbols $m$.

We conclude this section by showing that there are Delone sets
without uniform patch frequencies whose repetitivity functions
grow only slightly faster than linearly.
In particular, we show that Theorem~\ref{th43} is not far from best possible.

\begin{Lemma}\label{Le51}
Given positive integers $n\/$ and $N$, there exists a sequence
$\{a_k\}$ of positive integers with the following properties:
\begin{enumroman}
\item $a_k > N$ for all $k$,
\item $k\log k(\log \log k)^2<a_k<2k\log k(\log \log k)^2$
for all sufficiently large $k$;
\item $a_k$ is an even $n$th power;
\item the sequence $\{\rho_k\}$ of real numbers in $[0,1]$
defined recursively by
\begin{eqnarray}\label{eq604}
\rho_0&=&1,\nonumber\\
\rho_k&=&\frac{N}{2a_k}\rho_{k-1}+\left(1-\frac{N}{2a_k}\right)(1-\rho_{k-1})
\quad k=1,2,\ldots
\end{eqnarray}
does not tend to a limit.
\end{enumroman}
\end{Lemma}

\paragraph{Proof.} The recurrence (\ref{eq604}) can be rewritten as
$$
\rho_k-\frac{1}{2}=
-\left(1-\frac{N}{a_k}\right)\left(\rho_{k-1}-\frac{1}{2}\right),
$$
whence
$$
\rho_k=\frac{1}{2}+\frac{(-1)^k}{2}\prod_{j=1}^k\left(1-\frac{N}{a_j}\right).
$$
If we choose the sequence $\{a_k\}$ so that $a_k > N$ for each $k\/$ and
the infinite product
$$
\prod_{k=1}^\In\left(1-\frac{N}{a_j}\right)
$$
does not diverge to zero, then $\liminf \rho_k < \frac{1}{2}$ and
$\limsup \rho_k > \frac{1}{2}$.
Any choice of $a_k$'s satisfying (i), (ii) and (iii) of the lemma
(which are clearly compatible) fulfills these requirements.
\hspace*{\fill}$\Box$\smallskip

\begin{Theorem}\label{th52}
The $n\/$-dimensional cubic lattice $\ZZ^n$ can be two-colored in such a
way that, as $T \to \In$ the proportion of white lattice points in the cube
with center {\bf 0} and side length $T\/$ does not tend to a limit and
$$
M(T)=O(T(\log T)^{2/n}(\log\log\log T)^{4/n}).
$$
Here $M(T)$ is the natural extension to the two-colored lattice of the
repetitivity function for Delone sets of finite type.
\end{Theorem}

\paragraph{Proof.}  As a preliminary, let $\{a_k\}$ be a sequence
of positive integers satisfying Lemma~\ref{Le51} with $N = 2^{2^n+n}$,
and let $\{\rho_k\}$ be the corresponding sequence of real numbers.

It is convenient to regard the coloring as tiling $\RR^n$ with black
and white unit cubes with vertices at the points of the cubic lattice.

First take $2^{2^n}$ cubes of side 2, each composed of $2^n$ unit cubes and
each representing a different one of the $2^{2^n}$ possible vertex types at
the central vertex, and arrange them within the cube $\sC_1$ of side length
$s_1 = a_1^{1/n}$ that lies in the corner of the section of $\RR^n$ where
all coordinates are positive.  This can be done because $s_1$ is even and
$(s_1/2)^n > 2^{2^n}$, and it is possible to arrange that a white unit cube
is placed adjacent to the origin.  Fill the remainder of $\sC_1$ with black
unit cubes.  The proportion of white cubes in $\sC_1$ is $N/2a_1 = \rho_1$.
Next, repeat this construction with copies of $\sC_1$ in place of white
unit cubes, copies of the color-reversal $\overline{\sC}_1$ of $\sC_1$
in place of black unit cubes and $a_2$ in place of $a_1$.  The result
is that the cube $\sC_2$ of side length $s_2 = (a_1 a_2)^{1/n}$ at the
corner of the positive sector of $\RR^n$ is tiled with copies of $\sC_1$
and $\overline{\sC}_1$ in such a way as to include all `vertex types'
of them, to have a proportion $\rho_2$ of white unit cubes and to have
a copy of $\sC_1$ adjacent to the origin.  This construction can be
iterated indefinitely to give a tiling of the positive sector of $\RR^n$
by black and white unit cubes that, for each $k$, can be decomposed into
copies of $\sC_k$ and $\overline{\sC}_k$, where $\sC_k$ is the cube of
side length $s_k = (a_1 a_2 \cdots a_k)^{1/n}$ in the corner of the sector.
Finally, we extend the tiling to the whole of $\RR^n$ by reflection in
the coordinate hyperplanes.  The proportion of white unit cubes in
$\sC_k$ is $\rho_k$ and, since this does not tend to a limit, the
proportion of white unit cubes of the tiling in a large central cube
does not tend to a limit.

To bound $M(T)$ for the tiling, let $k\/$ be the smallest integer with
$s_k \ge 2T$.  Then, for any $\bx \in \ZZ^n$, $B(\bx,T)$ is contained in
a cube of side $2s_k$ made up of $2^n$ copies of cubes of types $\sC_k$
and $\overline{\sC}_k$ with a common vertex.  Since this vertex type occurs
in every copy of $\sC_{k+1}$ and $\overline{\sC}_{k+1}$ at the next level,
there is a point with the same $T$-patch as $\bx$ within every ball of
radius $\sqrt n\,s_{k+1}$.  By our choice of $k$, $s_{k-1} < 2T$, so
\begin{eqnarray}\label{eq509}
M(T)\le\sqrt n\,s_{k+1}&=&\sqrt n\,(a_k a_{k+1})^{1/n}s_{k-1}\nonumber\\
&<&2\sqrt n\,(a_k a_{k+1})^{1/n}T\nonumber\\
&=&O(T(k\log k(\log\log k)^2)^{2/n}),
\end{eqnarray}
\sloppy{by the right-hand inequality of Lemma~\ref{Le51}(ii).
Also, by the left-hand inequality of Lemma~\ref{Le51}(ii),}
$$
\log s_{k-1}>\frac{1}{n}\log((k-1)!)-O(1)>\frac{1}{n}k\log k-O(k)
$$
and hence
$$
k = O(\log T/\!\log\log T).
$$
So (\ref{eq509}) gives
$$
\hspace{40mm}M(T)=O(T(\log T)^{2/n}(\log\log\log T)^{4/n}).\hspace{37mm}\Box
$$
\begin{Corollary}\label{Co52}
There exists a Meyer set $X$ in $\RR^n$ with
$$
M_X (T)=O(T(\log T)^{2/n}(\log\log\log T)^{4/n})
$$
that does not have uniform patch frequencies.
\end{Corollary}

\paragraph{Proof.}  Replace each black lattice point \bx, in the coloring
given by Theorem~\ref{th52}, by the pair of points $\bx\pm(1,1,\ldots,1)/3$,
leaving the white lattice points unchanged.  The result is a Meyer set $X$ in
which the original white points can be identified by their $\sqrt n$-patches.
An argument similar to that at the end of the proof of Theorem~\ref{th51}
shows that if this patch had a uniform frequency $f$ in $X$ then the
frequency of its occurrence in the cube with center $\bf 0$ and side length
$T$ would tend to $f$ as $T\to\In$.  Since this is not the case, $X$ does
not have uniform patch frequencies.  Also $M_X(T)\le M(T+\sqrt n/3)$, where
the function on the right is the repetitivity function of the two-colored
lattice of Theorem~\ref{th52}.\hspace*{\fill}$\Box$\smallskip

\paragraph{Remarks.} The following three remarks clarify some features
of Theorem~\ref{th52} and its proof.

1. There are innumerable ways of coding the two colors by placement of points.
Perhaps a more elegant one is simply to translate all the black points by
$(1,1,\ldots,1)/2$, but the justification of this requires consideration
of $s_1\sqrt n$-patches, instead of $\sqrt n$-patches, because of the
impossibility of distinguishing between black and white points at smaller
scales.

2. The exponent of $\log T$ in Theorem~\ref{th52} is $2/n$ instead of
$1/n$ because rounding effects oblige us essentially to deal with cubes at
levels differing by two intead of at adjacent levels only.  One might think,
therefore, that the result of Corollary~\ref{Co52} would hold with an exponent
$1/n$ in place of $2/n$.  This is not so, however, because Theorem~\ref{th43}
can easily be strengthened to show that every Delone set $X\/$ with
$M_X(T)=O(T(\log T\log\log T )^{1/n})$ has a uniform average for every
local weight distribution.  The exponent $4/n\/$ of $\log\log\log T$
can, of course, be reduced to $(2+\eps)/n\/$ for any $\eps > 0$.

3. The upper and lower frequencies of white lattice points in cubes
with center $\bf 0$ are not local isomorphism invariants.
For example, if in the proof of Theorem~\ref{th52} we extend our
tiling of the positive sector using color-reversing reflections in
the coordinate hyperplanes, then we obtain a coloring of $\ZZ^n$
which is locally isomorphic to the coloring in the theorem but has
a proportion $1/2$ of white tiles in all boxes centered at the origin.

\section{Repetitive sets and the address map}

The address map introduced in \cite{Lag} provides a measure of the 
regularity of structure of finitely generated Delone sets $X$.

\paragraph{Definition 7.1.}
Let $X$ be a finitely generated Delone set in $\RR^n$, and suppose that the
additive group $[X]$ generated by $X$ in $\RR^n$ is of finite rank $s \ge n$.
An {\em address map} $\phi:[X]\to\ZZ^s$ is an isomorphism obtained 
by picking a fixed basis $\{\by_1,\ldots,\by_s\}$ of $[X]$.
Every $\by \in [X]$ has a unique decomposition
$\by =n_1\by_1+\cdots+n_s\by_s$, with all $n_i \in \ZZ$, and
\beql{eq702}
\phi(\by):=(n_1,n_2,\ldots,n_s).
\eeq
is the address map associated to this basis.

Address maps with different bases differ by left multiplication
by an element of $GL(s,\ZZ)$.
An address map is linear on $[X]$, which typically is dense in $\RR^n$,
but does not generally extend to a linear function on $\RR^n$.

In this section we show that the address maps of densely repetitive
and linearly repetitive sets satisfy (\ref{eq112}) and (\ref{JL113}),
respectively.

\begin{Theorem}\label{th71}
Let $X$ be a Delone set of finite type with ${\bf 0} \in X$ and
$\rank(X)=\rank(X-X)=s$.  Fix an address map $\phi:[X]\to\ZZ^s$.
If $X$ is either linearly repetitive or densely repetitive, then
there is a linear function $L:\RR^n\to\RR^s$ such that
\beql{eq706}
\|\phi(\bx)-L(\bx)\|=o(\|\bx\|)\quad\mbox{for $\bx\in X$}.
\eeq
In the case that $X$ is linearly repetitive, the error term
can be improved to $O(\|\bx\|^{1-\de})$, for some $\de=\de(X)>0$.
\end{Theorem}

\paragraph{Remark.}
The linear function $L\/$ of this theorem is uniquely determined
by $X$, since if $L^*$ is another such linear function then 
\beql{eq707}
\|L(\bx)-L^*(\bx)\|=o(\|\bx\|)\quad\mbox{for $\bx\in X$},
\eeq
giving $L\equiv L^*$ since $X\/$ is relatively dense.

\paragraph{Proof.}
We derive $L\/$ from the address map $\phi$ by a limiting process.
Given $\bt\in\RR^n$ let $\bx[\bt]$ denote some point in $X\/$ with
\beql{eq708}
\|\bx[\bt]-\bt\|\leq R.
\eeq
Since there are finitely many points of $X\/$ within a radius $R\/$
of $\bt$ it is clearly possible to adopt a convention for deciding,
in all circumstances, which one to pick, so resolving the ambiguity
of this definition.  Alternatively, we could note that the effects
of different choices are always absorbed by the error terms.

We shall define $L\/$ as a `smoothing' of $\phi\circ\bx$ by
\beql{eq709}
L(\bv):=\lim_{U\to\In}\frac{1}{U}\phi(\bx[U\bv])
\eeq
for $\bv\in\RR^n$.
In showing that this limit exists, it is convenient to prove a stronger
uniformly translation-invariant result, namely that the limit
\beql{eq710}
L(\bv)=\lim_{U\to\In}\frac{1}{U}\bp(\bt,U\bv+\bt)
\eeq
exists, for each $\bv\in\RR^n$, uniformly for $\bt\in\RR^n$
and is independent of $\bt$, where
$$
\bp(\bt_1,\bt_2)=\phi(\bx[\bt_2])-\phi(\bx[\bt_1]).
$$
represents the displacement in address space due to traversing a path
from $\bt_1$ to $\bt_2$ in object space.  (So (\ref{eq709}) is the
`centered' case $\bt = \bf 0$.)  It is sufficient to prove that the
limit (\ref{eq710}) exists when $\bv$ belongs to the standard basis
$\{\be_1,\ldots,\be_n\}$ of $\RR^n$, since if $\bv=\sum_{k=1}^n v_k\be_k$
and we put $\bv_i:=\sum_{k=i}^n v_k\be_k$ then
\beql{eq711}
\bp(\bt,U\bv+\bt)=\sum_{i=1}^n\bp(U\bv_{i+1}+\bt,Uv_i\be_i+U\bv_{i+1}+\bt).
\eeq
Moreover, if we can show that
\beql{eq712}
\frac{1}{U}\bp(\bt,U\be_i+\bt)=L(\be_i)+O(U^{-\de})
\eeq
uniformly in $\bt$, then we obtain (\ref{eq706}) with error term
$O(\|\bx\|^{1-\de})$ by taking $\bv=\bx$, $\bt=\bf 0$ and $U=1$ in
(\ref{eq711}).

The proof of (\ref{eq710}) is in two stages.
In the first stage, we replace the single path by a bunch of parallel
`sample paths' filling a box and use Theorems~\ref{th41} and \ref{th42} to
prove the existence of an average proportional displacement in address space.
In the second stage we apply this to a narrow box and use the Lipschitz
property to show that the average displacement is a good approximation
to the actual displacement due to traversing a single path.

For the first stage, let
$$
B=\{(t_1,t_2,\ldots,t_n):a_k\le t_k\le b_k\;(k=1,\ldots,n)\}
=[a_i,b_i]\times B^{(i)}
$$
be a box in $\RR^n$, where $B^{(i)}$ is the $(n-1)$-dimensional box
$$
B^{(i)}=\{(t_1,\ldots t_{i-1},t_{i+1},\ldots,t_n):
a_k\le t_k\le b_k\;(k=1,\ldots, i-1,i+1,\ldots n)\},
$$
let $G^{(i)}=B^{(i)}\cap\ZZ^{n-1}$ be the grid of points of
the $(n-1)$-dimensional integer lattice that are in $B^{(i)}$
and let $m_i$,$n_i$ be the integer parts of $a_i$ and $b_i$.
We have
\begin{eqnarray}\label{eq713}
|G^{(i)}|&=&\vol_{n-1}(B^{(i)})+O(\sig(B^{(i)}))\nonumber\\[1mm]
&=&\frac{\vol_n(B)}{b_i-a_i}\left(1+O\left(\frac{1}{\om(B)}\right)\right),
\end{eqnarray}
where $|G^{(i)}|$ is the cardinality of $G^{(i)}$.  Define
\beql{eq714}
\bP_i(B)=\sum_{\bg\in G^{(i)}}\bp((m_i,\bg),(n_i,\bg))
\eeq
with the convention that when $\bg$ is on the boundary of $B^{(i)}$, the
term for $\bg$ is to be taken with an appropriate weight---$\frac{1}{2}$
when $\bg$ is on an $(n-2)$-dimensional face, $\frac{1}{4}$ when it is on
an $(n-3)$-dimensional face, and so on.
The sum $\bP_i(B)$ represents the sum of the displacements in address space
due to traversing a number of paths in object space parallel to $\be_i$.
The endpoints of these paths are points of the integer lattice $\ZZ^n$
close to the faces of $B\/$ that are perpendicular to $\be_i$.
It is readily seen that $\bP_i$ is a local weight distribution on $X$.
Definition~5.2(a) is a consequence of the Lipschitz property (\ref{eq201})
of $\phi$ and our estimate for $|G^{(i)}|$.  Part~(c) is satisfied without
the error term (when the weights assigned to terms in the sum for $\bP_i(B)$
arising from points on the boundary of $B^{(i)}$ are taken into account).
The error in part~(b) is due to the fact that, although $\bx[(g,\bg)]$
lies within a radius $R\/$ of $(g,\bg)$, which point of $X\/$ it is may
depend on nearby points of $X\/$ outside $B$.
In view of the Lipschitz condition, the total error is
$$
O(|G^{(i)}|)=O(\vol_{n-1}(B^{(i)}))=O(\sig(B)).
$$
Now Theorems \ref{th41} and \ref{th42} ensure the existence
of a limiting value $L(\be_i)$ such that
\beql{eq715}
\frac{\bP_i(B)}{\vol_n(B)}=L(\be_i)+\eps(B)
\eeq
with $\eps(B)\to0$ as $\om(B)\to\In$ in both the linearly and
densely repetitive cases.  In the linearly repetitive case
$\eps(B)=O(\om(B)^{-\de_1})$, where $\de_1$ is the $\de$ of
Theorem~\ref{th41}.

The second stage of the proof is to use this average value of
$\bP_i$ to estimate $\bp(\bt,U\be_i+\bt)$ for a sample path
by applying it to the box $B\/$ defined by
$$
[a_k,b_k] = \cases{[t_i,t_i+U]& for $k=i$,\vspace{2mm}\cr
[t_k-\frac{1}{2}\sqrt U,t_k+\frac{1}{2}\sqrt U]& for $k\ne i$,\cr}
$$
where $t_1\be_1+\cdots+t_n\be_n=\bt$.  By the Lipschitz property
(\ref{eq201}), every term on the right-hand side of (\ref{eq714}) is
$$
\bp(\bt,U\be_i+\bt)+O(\sqrt U).
$$
So (\ref{eq714}), (\ref{eq713}) and the Lipschitz property give
$$
\bP_i(B)=U^{(n-1)/2}\bp(\bt,U\be_i+\bt)+O(U^{n/2})
$$
and by (\ref{eq715})
$$
\frac{\bp(\bt,U\be_i+\bt)}{U}=L(\be_i)+\eps(B)+O(U^{-1/2}),
$$
with $\eps(B)\to0$ as $U\to\In$ in both the linearly and densely repetitive
cases and $\eps(B) = O(U^{-\de_1/2})$ in the linearly repetitive case.
This establishes the existence of the limit (\ref{eq710}) in both cases
and gives the estimate (\ref{eq712}) in the linearly repetitive case
with $\de=\de_1/2$.\hspace*{\fill}$\Box$\smallskip

The following example shows that for linearly repetitive sets $X$ the bound
$$\|\phi(\bx)-L(\bx)\|=O(\|\bx\|^{1-\de})$$
given in Theorem~\ref{th71} cannot be made smaller than a power of $\|\bx\|$.

\begin{Example}
There is a linearly repetitive Delone set $X$ in $\RR^n$ which has
the property that there is a constant $\de=\de(X)<1$ and a sequence
$\{\bx_i\}\subseteq X$ with $\|\bx_i\|\to\In$, such that
\beql{eq634}
\|\phi(\bx_i)-L(\bx_i)\|>\|\bx_i\|^{1-\de}.
\eeq
\end{Example}

\paragraph{Proof.}
Such a set $X$ is constructed as a set of control points of a
self-similar tiling that has an expansion constant $\beta$ which
is not a Pisot or Salem number, see Kenyon~\cite{Ken96}.  The set
of control points of any self-similar tiling is a linearly repetitive
set, according to a result of Solomyak \cite[Lemma~2.3]{Sol98a}.
The fact that the expansion constant $\beta$ has an algebraic conjugate
with absolute value greater than one produces property (\ref{eq634}). 
It follows that such a set $X$ cannot be a Meyer set, because (\ref{eq634})
violates the Meyer set condition (\ref{eq202}).\hspace*{\fill}$\Box$

\section{Linear repetitivity, diffraction and model sets}

In this section we study linearly repetitive Delone sets, and
consider them as models for `perfectly ordered quasicrystals'.

We begin by stating a theorem proved in \cite{LP}, which shows that
linear repetitivity is the slowest possible growth rate of $M_X(T)$
for any translation-minimal set $X$ that is not fully periodic. 

\begin{Theorem}\label{th81} 
Let $X$ be a repetitive Delone set in $\RR^n$.
If its repetitivity function $M_X(T)$ satisfies
\beql{801a}
M_X(T) < \frac {1}{3} T
\eeq
for a single value of $T > 0 $, then $X$ is an ideal crystal.
\end{Theorem}

This result implies that Conjectures~1.2a and 1.2b are equivalent, as follows.
Any linearly repetitive set $X$ in $\RR^n$ has $N_X(T) = O (T^n)$ by
Theorem~\ref{th31}.  Conjecture~1.2a asserts that if $X$ is aperiodic
then it is densely repetitive, which gives  
\beql{801b}
M_X(T)^n < c N_X(T) \quad \mbox{for $T>T_0$},
\eeq
with some positive $c$.
Now Theorem~\ref{th81} gives the linear lower bound $M_X(T)>T/3$ and
(\ref{801b}) yields $N_X(T) > c'T^n$ for $T>T_0$, which is Conjecture~1.2b.
In the other direction, if Conjecture~1.2b holds, then an aperiodic
linearly repetitive set has $N_X(T) > c' T^n$ for $T>T_0$ and linear
repetitivity then gives $M_X(T) < c_1 T < c_2 (N_X(T))^{1/n}$, so $X$
is densely repetitive, which is Conjecture~1.2a.

We next establish some existence results for linearly repetitive sets.
We construct a large class of such sets in the one-dimensional case, 
given in Theorem~\ref{th82} below.  It seems much harder to construct
linearly repetitive sets in dimensions $n \geq 2$ that are irreducible in the
sense of not being a product of lower-dimensional linearly repetitive sets.
There do exist non-trivial examples of linearly repetitive sets in
all dimensions $n \geq 2$, based on self-similar constructions. 

The one-dimensional construction is based on Beatty sequences of numbers 
$\alpha \in [0,1]$ whose continued fractions have bounded partial quotients.
Given any symbol sequence $S=\{s_i\}\in\{0,1\}^\ZZ$ and any $\tau > 1$,
we obtain a set $X_\tau (S)$ by assigning $x_0 = 0$ and then measuring
distances $x_{i+1}-x_i$ to be 1 or $\tau$ according as $s_i$ is 0 or 1.
The set $X_\tau (S)$ is a Delone set of finite type and it is easy to
see that whether or not $X_\tau (S)$ is linearly repetitive depends
only on the symbol sequence $S$.

\begin{Lemma}\label{le81}
Let $S \in \{0,1\}^\ZZ$.
Then for any $\tau>1$, $X_\tau (S)$ is linearly repetitive if and
only if $S$ is minimal with recurrence function $\tilde{M}_S (n)$
satisfying $\tilde{M}_S (n) = O(n)$ as $n\to\In$.
\end{Lemma}

\paragraph{Proof.}
This is immediate.\hspace*{\fill}$\Box$\smallskip

The aperiodic 2-symbol sequences $S$ with the smallest word-counts have
exactly $n+1 $ words of length $n$ for each $n\geq1$, and of these, the
ones which are minimal are exactly the Beatty sequences of irrational slope,
see Coven and Hedlund~\cite{CoHe73}, Coven~\cite{Cov75} and Paul~\cite{Pau75}.
It is easy to characterize those Beatty sequences that have a recurrence
function with linear growth $O(n)$.

\begin{Theorem}\label{th82} 
Let $B_\alpha \subseteq \{0,1\}^\ZZ$ be a Beatty sequence with
$\alpha\in[0,1]$.  Then:
\begin{enumroman}
\item the recurrence function $\tilde{M}_{B_\alpha}(T)$
is bounded if and only if $\alpha$ is rational;
\item the recurrence function satisfies
$$\tilde{M}_{B_\alpha}(T) = O(T)\quad\mbox{as }T\to\In$$ 
if and only if $\alpha$ is badly approximable, i.e.\ $\alpha$ is a number
whose regular continued fraction expansion has bounded partial quotients.
\end{enumroman}
\end{Theorem}

Recall that badly approximable numbers are characterized as the
numbers $\alpha$ for which there exists a constant $c>0$ such
that $|\alpha-a/q|>c/q^2$ for every rational number $a/q$.
For their various properties see Shallit~\cite{Sha92}.

\paragraph{Proof.}
When $\alpha$ is rational with denominator $q$, then $B_\alpha$ has
period $q\/$ and $\tilde{M}_{B_\alpha}(T) = q\/$ for all $T\ge q$.

Now let $a_k$ be the $k\/$th partial quotient and $p_k/q_k$ be the $k\/$th
convergent of the regular continued fraction expansion of $\alpha$.
By (\ref{eq407}) we have
$$
\tilde{M}_{B_\alpha}(\ell)=q_k+q_{k+1}\le (a_{k+1}+2)q_k\le(a_{k+1}+2)\ell
$$
for $q_k\leq\ell<q_{k+1}$ and
$$
\tilde{M}_{B_\alpha}(q_k)=q_k+q_{k+1}\ge(a_{k+1}+1)q_k.
$$
Also $q_k\to\In$ as $k\to\In$ when $\alpha$ is irrational and,
in particular, when $\{a_k\}$ is unbounded.
Hence, when $\alpha$ is irrational $\tilde{M}_{B_\alpha}(T)=O(T)$
is equivalent to the sequence of partial quotients $\{a_k\}$ of
$\alpha$ being bounded.\hspace*{\fill}$\Box$\smallskip

Theorem~\ref{th82} implies that there are restrictions on the image in
$\RR^s$ of the linear form $L$ that is associated by Theorem~\ref{th71}
to a linearly repetitive Delone set $X_{\tau}(B_{\alpha})$: this image is
necessarily a line in $\RR^2$ whose slope is a badly approximable number.

The aperiodic Delone sets $X$ in $\RR^n$ that are linearly repetitive
have several claims to be the `simplest' aperiodic sets.
These claims include the following.
\begin{enumarabic}
\item Their patch-counting functions satisfy $N_X (T) = O(T^n)$.
If Conjecture~\ref{Cj21} is correct, this is the slowest possible
growth rate of $N_X (T)$ for aperiodic sets.
\item Among all aperiodic sets, they have minimal growth rate of
the repetitivity function $M_X (T)$.
\item They have long-range order in the sense of having strict
uniform patch frequencies and, hence, are diffractive.
\item They include most, if not all, of the proposed examples of strongly
ordered sets with quasicrystalline properties, including point sets
associated to irreducible self-affine tilings \cite[Lemma~2.3]{Sol98a}.
They also include various special cut-and-project sets such as Penrose
tilings.
\end{enumarabic}

For these reasons we consider the class of linearly repetitive Delone sets
to be a candidate for the concept of {\em perfectly ordered quasicrystals}.
This class includes ideal crystals; we do not impose aperiodicity as
a requirement for being quasicrystalline.

Many authors, however, consider that any definition of quasicrystalline
order should also require a well-defined diffraction measure with a
non-trivial pure point component (`Bragg peaks').

\paragraph{Definition~8.1}
A Delone set $X$ is said to have {\em strong long range order} or to
be {\em strongly diffractive} if it has a unique diffraction measure
$\hat{\gamma}_X$ which includes a pure point component that is
relatively dense in $\RR^n$.

Not all linearly repetitive Delone sets have this property.
The construction of Example~7.1 can be used to produce a linearly
repetitive Delone set $X$ which is not a Meyer set and which has no
non-trivial pure point spectrum, using the criterion of Solomyak
\cite[Theorem~2.1]{Sol98b}, so is not strongly diffractive.
In the opposite direction, cut-and-project sets with a nice window
are always strongly diffractive, but need not be linearly repetitive.

We may further clarify the relations between linear repetitivity, 
strong diffraction, and being a Meyer set, by considering some
properties of the address map and the associated topological
dynamical system of a finitely generated Delone set $X$. 
The orbits of the topological dynamical system $(\sX_X,\RR^n)$ under the
$\RR^n$-action can be described using the address map $\phi:[X-X]\to\ZZ^s$
and a projection map $\pi:\RR^s\to\RR^n$ (not necessarily orthogonal).
We have the following commutative diagram:
\beql{XX701}
\begin{array}{c}
\setlength{\unitlength}{1mm}
\begin{picture}(40,20)
\put(10,2){\makebox(0,0){$[X-X]$}}
\put(10,18){\makebox(0,0){$\ZZ^s$}}
\put(30,2){\makebox(0,0){$\RR^n$}}
\put(30,18){\makebox(0,0){$\RR^s$}}
\put(6,10){\makebox(0,0){$\phi$}}
\put(32,10){\makebox(0,0){$\pi$}}
\put(20,2){\vector(1,0){6}}
\put(20,2.8){\oval(2,1.6)[l]}
\put(16,18){\vector(1,0){10}}
\put(16,18.8){\oval(2,1.6)[l]}
\put(9.5,5){\vector(0,1){10}}
\put(29,15){\vector(0,-1){10}}
\end{picture}
\end{array}
\eeq
Here the horizontal arrows are inclusions and $\pi$ is the projection
map, defined by $\pi(\be_j)=\bv_j$ for $1\le j\le s$ where
$\{\bv_1,\ldots,\bv_s\}$ is the chosen basis of $X-X$ and $\be_j$
are the unit coordinate vectors. If $L:\RR^n\to\RR^s$ is the linear
map given by Theorem~\ref{th71}, then $\pi\circ L$ is the identity.
If ${\bf 0}\in X$, then the addresses $Y=\phi (X)$ of $X$ together
with $\pi$ permit $X$ to be uniquely reconstructed using (\ref{XX701}).
If $(\sX_X,\RR^n)$ is minimal, then $[X'-X']$ is independent of the
choice of $X'\in\sX_X$ and a general $X'\in\sX_X$ is described by
the data $(Y',\pi,\bx')$, where we choose any point $\bx' \in X'$
and $Y'=\phi(X'-\bx')$ has ${\bf 0}\in Y'\subseteq\phi[X'-X']$.
This is unique up to the choice of $\bx'$, so the $\RR^n$-orbit of
$X'$ is specified by the equivalence class $\{ Y'\}$ of translates
of $Y'$ in $\ZZ^s$ that contain $\bf 0$, i.e.
$$
\{Y'\}:=\{Y'-\bm:\bm\in Y'\}.
$$
We regard $\{Y'\}$ as describing a `symbolic dynamics' for this orbit.

The properties of minimality and linear repetitivity of
$(\sX_X,\RR^n)$ are completely determined by the equivalence
class $\{Y\}$ of $Y:=\phi(X)$ (where we assume ${\bf 0}\in X\/$).
We can now obtain a new topological dynamical system by holding
$Y$ fixed and varying the projection $\pi$, to $\pi'$ say.
The pair $(Y,\pi')$ uniquely determines the preimage $X'$ by
(\ref{XX701}), and if $X'$ is a Delone set  then we obtain
a new Delone dynamical system $(\sX_{X'},\RR^n)$.
If $Y$ is linearly repetitive then the new topological dynamical
system $(\sX_{X'},\RR^n)$ is uniquely ergodic, hence gives rise
to a unique metric dynamical system $(\sX_{X'},\RR^n,d\mu')$,
and all sets in $\sX_{X'}$ are diffractive.

On the other hand, the existence of a pure point component of the
diffraction measure $\hat{\gamma}_{X'}$ depends on the projection $\pi'$.
Kol\'{a}\v{r} {\em et~al} \cite[\S5.2]{Kol93} have one-dimensional
examples in which $Y \subseteq \ZZ^2$ is linearly repetitive and
$X'=\pi(Y)$ seems to have some pure point spectrum for certain
projections $\pi$, but no pure point spectrum for other projections.
(In this example the sets $X'$ correspond to tilings of the line with
tiles of two lengths and varying $\pi$ is equivalent to varying the ratios
of the tile lengths while leaving the  pattern of the tiles fixed.)

To summarize: {\em the property of being linearly repetitive depends
only on the symbolic equivalence class $\{Y\}$, while the presence or
absence of a pure point component in the diffraction spectrum depends
on both $\{Y\}$ and the projection $\pi$.}

We now look at some one-dimensional examples, given by sets $X_\tau(S)$
constructed, as before, from a symbol sequence $S$ and a real parameter
$\tau$. 
Lemma~\ref{le81} showed that the linear repetitivity of $X_\tau(S)$
is determined by $S$ alone and is independent of $\tau$.
In contrast, the property of $X_\tau (S)$ being a Meyer set sometimes depends
on whether $\tau$ is rational or irrational, see \cite[Section~5]{Lag}.
This can occur because when $\tau$ is rational $\phi$ maps $[X]$
to $\ZZ$, but when $\tau$ is irrational it maps $[X]$ to $\ZZ^2$.
For Beatty sequences $B_{\alpha}$, the sets $X_{\tau}(B_{\alpha})$
are always Meyer sets by the criterion \cite[Theorem 5.1]{Lag}.

We end with some open problems suggested
by the observations made in this section.

\paragraph{Problem~8.1.}
Explicitly characterize all two-sided shifts $\Sigma$ on a finite alphabet
that are minimal and have recurrence function $\tilde{M}_{\Sigma}(\ell)$
bounded above by $C\ell\/$ for some constant $C$.

\paragraph{Problem~8.2.}
Characterize the possible images in $\RR^s$ of linear forms
$L:\RR^n\to\RR^s$ for which there exists a linearly repetitive
Delone set $X$ with ${\bf 0}\in X$ that is associated to $L\/$
in the manner of Theorem~\ref{th71}.

\paragraph{Problem~8.3}
Explicitly characterize all cut-and-project sets that are linearly repetitive.
\bigskip

\paragraph{Acknowledgements.}
We are indebted to M.~Baake, B.~Solomyak and the referee for helpful
comments and to G.~van Ophuysen for pointing out results of Morse and
Hedlund relevant to Theorems~\ref{th32} and~\ref{th82}.

\section*{{\rm A.}\quad Appendix. Repetitivity function for cubical patches}

We can define the {\em cubical repetitivity function} $\tilde{M}_X(T)$
to be the smallest value of $M\/$ such that every cube of side $2M+T\/$
contains a translate of every cubical $T$-patch of $X$.
Let $\tilde{N}_X(T)$ count the number of cubical patches of side $T$.
Then a Delone set $X$ is a {\em densely repetitive set for cubes\/} if
$$
\tilde{M}_X(T) = O(\tilde{N}_X(T)^{1/n}),\quad\mbox{as $T\to\In$.}
$$
This is equivalent to $X$ being densely repetitive if $N_X(T)$ grows
polynomially in $T$, or if the dimension $n =1 $; for $N_X(T)$ and
$\tilde{N}_X(T)$ are then within a constant factor of each other.
However, if $N_X (T)$ grows faster than polynomially the two concepts
may differ, as far as we know.
The proof of Theorem~\ref{th51} also works in the context of densely
repetitive Delone sets for cubes and shows that they too have strict
uniform patch frequencies.

\catcode`\@=11 
\renewcommand{\section}{%
\@startsection{section}{1}{0in}{-.5in}{.25in}{\centering\sc}}
\catcode`@=12


\begin{thebibliography}{\bf 99}

\footnotesize

\bibitem{All94}
J.-P.~Allouche.
Sur la complexit\'e des suites infinies.
{\em Bull.\ Belg.\ Math.\ Soc.\ Simon Stevin} {\bf 1} (1994), 133--148.

\vspace*{-3mm}
\bibitem{Be92}
J.~Bellissard.
Gap Labelling Theorems for Schr\"odinger Operators.
{\em From Number Theory to Physics}.
Eds.\ M.~Waldschmidt, P.~Moussa, J.-M.~Luck and C.~Itzykson.
Springer, Berlin, 1992, pp.~538--630.

\vspace*{-3mm}
\bibitem{BV99}
V.~Berth\'e and L.~Vuillon.
Suites doubles de basse complexit\'e.
{\em J.~Th\'{e}or.\ Nombres Bordeaux} {\bf 12} (2000), 179--208.

\vspace*{-3mm}
\bibitem{Bos84}
M.~Boshernitzan,
A unique ergodicity of minimal symbolic flows with linear block growth,
{\em J.~Anal.\ Math.}\ {\bf 44} (1984/85), 77--96.

\vspace*{-3mm}
\bibitem{Cov75}
E.~M.~Coven.
Sequences with minimal block growth.
{\em Math.\ Syst.\ Th.}\ {\bf 8} (1975), 376--382.

\vspace*{-3mm}
\bibitem{CoHe73}
E.~M.~Coven and G.~A.~Hedlund.
Sequences with minimal block growth.
{\em Math.\ Syst.\ Th.}\ {\bf 7} (1973), 138--153.

\vspace*{-3mm}
\bibitem{De}
B.~N.~Delone [B.~N.~Delaunay], N.~P.~Dolbilin,
M.~I.~Shtogrin and R.~V.~Galiulin.
A local criterion for regularity of a system of points.
{\em Sov.\ Math.\ Dokl.}\ {\bf 17}(2) (1976), 319--322.

\vspace*{-3mm}
\bibitem{DolLagSen}
N.~P.~Dolbilin, J.~C.~Lagarias and M.~Senechal.
Multiregular point systems.
{\em Discrete Comput.\ Geom.}\ {\bf 20} (1998), 477--498.

\vspace*{-3mm}
\bibitem{DolPle}
N.~P.~Dolbilin and P.~A.~B.~Pleasants.
Aperiodic sets with few isometry patches.
in preparation.

\vspace*{-3mm}
\bibitem{Els86}
V.~Elser.
The diffraction pattern of projected structures.
{\em Acta Cryst.\ Sect.\ A} {\bf 42} (1986) 36--43.

\vspace*{-3mm}
\bibitem{Fe96}
S.~Ferenczi.
Rank and symbolic complexity.
{\em Ergod.\ Th.\ \& Dynam.\ Sys.}\ {\bf 16} (1996), 663--682.

\vspace*{-3mm}
\bibitem{Fur67}
H.~Furstenberg.
Disjointness in ergodic theory, minimal sets,
and a problem in Diophantine approximation.
{\em Math.\ Syst\. Th.}\ {\bf 1} (1967), 1--49.

\vspace*{-3mm}
\bibitem{Fur78}
H.~Furstenberg.
{\em Recurrence in Ergodic Theory and Combinatorial Number Theory.}
Princeton University Press, Princeton, NJ, 1978.

\vspace*{-3mm}
\bibitem{GahKli}
F.~G\"ahler and R.~Klitzing.
The diffraction pattern of self-similar tilings.
{\em The Mathematics of Long-Range Aperiodic Order}.
Ed.\ R.~V.~Moody.  Kluwer, Dordrecht, 1997, pp.~141--174.

\vspace*{-3mm}
\bibitem{Getal}
A.~I.~Goldman {\em et al}.  
Quasicrystalline Materials.
Amer.\ Scient.\ {\bf 84} (1996), 230--241.

\vspace*{-3mm}
\bibitem{GruShe87}
B.~Gr\"unbaum and G.~C.~Shephard.
{\em Tilings and Patterns}.
W.~H.~Freeman \& Co., New York 1987.

\vspace*{-3mm}
\bibitem{He91}
C.~L.~Henley.
Random Tiling Models.
{\em Quasicrystals: The State of the Art}
Eds.\ D.~P.~DiVincenzo and P.~J.~Steinhardt.
World Scientific, Singapore, 1991, pp.~429--524.

\vspace*{-3mm}
\bibitem{Hof95}
A.~Hof.
On diffraction by aperiodic structures.
Commun.\ Math.\ Phys.\ {\bf 169} (1995), 25--43.

\vspace*{-3mm}
\bibitem{Hof95a}
A.~Hof.
Diffraction by aperiodic structures.
{\em The Mathematics of Long-Range Aperiodic Order}.
Ed.\ R.~V.~Moody.
Kluwer, Dordrecht, 1997, pp.~239--268.

\vspace*{-3mm}
\bibitem{Jan92}
C.~Janot.
{\em Quasicrystals: a Primer}.
Clarendon Press, Oxford 1992 (2nd edn 1994).

\vspace*{-3mm}
\bibitem{Ken92}
R.~Kenyon.
Self-replicating tilings.
{\em Symbolic Dynamics and its Applications}.
Ed.\ P.~Walters. 
American Mathematical Society, Providence, RI, 1992, pp.~239--264.

\vspace*{-3mm}
\bibitem{Ken96}
R.~Kenyon.
The construction of self-similar tilings.
{\em Geom.\ Funct.\ Anal.}\ {\bf 6} (1996), 471--488.

\vspace*{-3mm}
\bibitem{Kol93}
M.~Kol\`ar, B.~Iochum and L.~Raymond.
Structure factor of ID systems (superlattices) based on two-letter
substitution rules: I.~Bragg peaks.
{\em J.~Phys.\ A} {\bf 26} (1993), 7343--7366.

\vspace*{-3mm}
\bibitem{Lag96}
J.~C.~Lagarias. 
Meyer's concept of quasicrystal and quasiregular sets.
{\em Commun.\ Math.\ Phys.}\ {\bf 179} (1996), 365--376.

\vspace*{-3mm}
\bibitem{Lag}
J.~C.~Lagarias.
Geometric models for quasicrystals I.  Delone sets of finite type. 
{\em Discrete Comput.\ Geom.}\ {\bf 21} (1999), 161--191.

\vspace*{-3mm}
\bibitem{Lag98}
J.~C.~Lagarias.
Geometric models for quasicrystals II.  Local rules under isometries. 
{\em Discrete Comput.\ Geom.}\ {\bf 21} (1999), 345--372.

\vspace*{-3mm}
\bibitem{LP}
J.~C.~Lagarias and P.~A.~B.~Pleasants.
Local complexity of Delone sets and crystallinity.
{\em Canad.\ Math.\ Bull.}\ {\bf 45} (2002), 634--652.

\vspace*{-3mm}
\bibitem{LePiuSad93}
T.~Q.~T.~Le, S.~Piunikhin and V.~Sadov.
Geometry of quasicrystals.
{\em Russ.\ Math.\ Surveys} {\bf 48} (1993), 41--102.

\vspace*{-3mm}
\bibitem{LMS02}
J.-Y.~Lee, R.~V.~Moody and B.~Solomyak.
Pure point dynamical and diffractive spectra.
{\em Ann.\ Inst.\ H.~Poincar\'{e}} {\bf 3} (2002), 1003--1018.

\vspace*{-3mm}
\bibitem{Lenz02}
D.~Lenz.
Aperiodic linearly repetitive Delone sets are densely repetitive.
{\em Discrete Comput.\ Geom.}\ to appear.

\vspace*{-3mm}
\bibitem{LS02}
D.~Lenz and P.~Stollmann.
Delone dynamical systems and associated random operators.
{\em Proceedings: Operator Algebras and Mathematical Physics,
Constanta 2001.}\\
Eprint arXiv math-ph/0202042.

\vspace*{-3mm}
\bibitem{Lin84}
D.~A.~Lind.
The entropies of topological Markov shifts
and a related class of algebraic integers.
{\em Ergod.\ Th.\ Dynam.\ Syst.}\ {\bf 4} (1984), 283--300.

\vspace*{-3mm}
\bibitem{Lin92}
D.~A.~Lind.
Matrices of Perron numbers.
{\em J.~Number Th.}\ {\bf 40} (1992), 211-217.

\vspace*{-3mm}
\bibitem{LinMar95}
D.~A.~Lind and B.~Marcus.
{\em An Introduction to Symbolic Dynamics and Coding.}
Cambridge University Press, Cambridge, 1995.

\vspace*{-3mm}
\bibitem{LP1}
W.~F.~Lunnon and P.~A.~B.~Pleasants.
Quasicrystallographic tilings.
{\em J.~Math.\ Pures Appl.}\ {\bf 66} (1987), 217--263.

\vspace*{-3mm}
\bibitem{LP2}
W.~F.~Lunnon and P.~A.~B.~Pleasants.
Characterization of two-distance sequences.
{\em J.~Austral.\ Math.\ Soc.\ Ser.\ A} {\bf 53} (1992), 198--218.

\vspace*{-3mm}
\bibitem{Mey70}
Y.~Meyer.
{\em Nombres de Pisot, nombres de Salem, et analyse harmonique
(Lecture Notes in Mathematics, 117).} 
Springer, Berlin, 1970.

\vspace*{-3mm}
\bibitem{Mey72}
Y.~Meyer.
{\em Algebraic Numbers and Harmonic Analysis}.
North-Holland, Amsterdam, 1972.

\vspace*{-3mm}
\bibitem{Mey95}
Y.~Meyer.
Quasicrystals, Diophantine approximation and algebraic numbers.
{\em Beyond Quasicrystals}. 
Eds.\ F.~Axel and D.~Gratias.
Springer, Berlin, 1995, pp.~3--16.

\vspace*{-3mm}
\bibitem{Moo95}
R.~V.~Moody. 
Meyer sets and the finite generation of quasicrystals.
{\em Symmetries in Science VIII}.
Ed.\ P.~M.~Gruber. 
Plenum, London, 1995.

\vspace*{-3mm}
\bibitem{Moo97}
R.~V.~Moody. 
Meyer sets and their duals.
{\em The Mathematics of Long-Range Aperiodic Order}.
Ed.\ R.~V.~Moody.
Kluwer, Dordrecht, 1997, pp.~403--442.

\vspace*{-3mm}
\bibitem{MooPat95}
R.~V.~Moody and J.~Patera.
Colourings of quasicrystals.
{\em Canad.\ J.\ Phys.}\ {\bf 72} (1995), 442--452.

\vspace*{-3mm}
\bibitem{MorHed38}
M.~Morse and G.~A.~Hedlund.
Symbolic dynamics.
{\em Amer.\ J.\ Math.}\ {\bf 60} (1938), 815--866.

\vspace*{-3mm}
\bibitem{MorHed40}
M.~Morse and G.~A.~Hedlund.
Symbolic dynamics~II, Sturmian trajectories.
{\em Amer.\ J.\ Math.}\ {\bf 62} (1940), 1--42.

\vspace*{-3mm}
\bibitem{Moz89}
S.~Mozes.
Tilings, substitution systems and dynamical systems generated by them.
{\em J.~Anal.\ Math.}\ {\bf 53} (1989), 139--186.

\vspace*{-3mm}
\bibitem{OguDunKat88}
C.~Oguey, M.~Duneau and A.~Katz.
A geometric approach to quasiperiodic tilings.
{\em Commun.\ Math.\ Phys.}\ {\bf 118} (1988), 99--118.

\vspace*{-3mm}
\bibitem{Pa81}
W. Parry.
{\em Topics in Ergodic Theory.}
Cambridge University Press, Cambridge, 1981.

\vspace*{-3mm}
\bibitem{Pau75}
E.~M.~Paul.
Minimal symbolic flows having minimal block growth.
{\em Math.\ Syst.\ Th.}\ {\bf 8} (1975), 309--315.

\vspace*{-3mm}
\bibitem{Pet}
K.~Petersen.
{\em Ergodic Theory}.
Cambridge University Press, Cambridge, 1983.

\vspace*{-3mm}
\bibitem{Ple1}
P.~A.~B.~Pleasants.
Designer quasicrystals: cut-and-project sets with pre-assigned properties.
{\em Directions in Mathematical Quasicrystals (CRM Monograph Series).}
Eds.\ M.~Baake and R.~V.~Moody.
American Mathematical Society, Providence, RI, 2000.

\vspace*{-3mm}
\bibitem{Ple2}
P.~A.~B.~Pleasants.
Entropy of languages and sequences.
{\em Manuscript}, 1986.

\vspace*{-3mm}
\bibitem{Que87}
M.~Queffelec.
{\em Substitution Dynamical Systems---Spectral Analysis
(Lecture Notes in Mathematics, 1294)}
Springer, Berlin, 1987.

\vspace*{-3mm}
\bibitem{Que95}
M.~Queffelec.
Spectral study of automatic and substitutive sequences.
{\em Beyond Quasicrystals.}
Eds.\ F.~Axel and D.~Gratias.
Springer, Berlin, 1995, pp.~369--414.

\vspace*{-3mm}
\bibitem{Rad91a}
C.~Radin.
Disordered ground states for classical lattice models.
{\em Rev.\ Math.\ Phys.}\ {\bf 3} (1991), 125--135.

\vspace*{-3mm}
\bibitem{Rad91}
C.~Radin.
Global order from local sources.
{\em Bull.\ Amer.\ Math.\ Soc.}\ {\bf 25} (1991), 335--364.

\vspace*{-3mm}
\bibitem{RadWol92}
C.~Radin and M.~Wolff.
Space tilings and local isomorphism.
{\em Geom.\ Dedicata} {\bf 42} (1992), 355--360.

\vspace*{-3mm}
\bibitem{RHB98}
C.~Richard, M.~H\"offe, J.~Hermisson and M.~Baake.
Random tilings: concepts and examples.
{\em J.~Phys.\ A} {\bf 31} (1998), 6385--6408.

\vspace*{-3mm}
\bibitem{Rob9394}
E.~A.~Robinson, Jr.
The dynamical theory of tilings and quasicrystallography.
{\em Ergodic Theory of $\ZZ^d$-Actions 
(London Mathematical Society Lecture Notes, 228).} 
Eds. M.~Pollicott and K.~Schmidt.
Cambridge University Press, Cambridge, 1996, pp.~451--473.

\vspace*{-3mm}
\bibitem{Rob}
E.~A.~Robinson, Jr.
The dynamical properties of Penrose tilings.
{\em Trans.\ Amer.\ Math.\ Soc.}\ {\bf 348} (1996), 4447--4464.

\vspace*{-3mm}
\bibitem{ST99a}
J.~W.~Sander and R.~Tijdeman.
The complexity of functions on lattices.
{\em Th.\ Comput.\ Sci.}\ {\bf 246} (2000), 195--225.

\vspace*{-3mm}
\bibitem{ST99b}
J.~W.~Sander and R.~Tijdeman.
The rectangle complexity of functions on two-dimensional lattices.
{\em Th.\ Comput.\ Sci.}\ {\bf 270} (2002), 857--863.

\vspace*{-3mm}
\bibitem{Schlot}
M.~Schlottmann.
{\em Geometrische Eigenschaften quasiperiodischer Strukturen.}
Doctoral thesis.
Eberhard-Karls-Universit\"at zu T\"ubingen, Germany, 1993.

\vspace*{-3mm}
\bibitem{Sch}
M.~Schlottmann.
Cut and project sets in locally compact Abelian groups.
{\em Quasicrystals and Discrete Geometry (Fields Institute Monographs, 10).}
Ed.\ J.~Patera.
American Mathematical Society, Providence, RI, 1998, pp.~247--264.

\vspace*{-3mm}
\bibitem{Sch93}
J.-P.~Schreiber. 
Approximations diophantiennes et probl\`ems additifs
dans les groupes ab\'eliens localement compacts.
{\em Bull.\ Soc.\ Math.\ France} {\bf 101} (1973), 297--332.

\vspace*{-3mm}
\bibitem{Sen95}
M.~Senechal.
{\em Quasicrystals and Geometry}.
Cambridge University Press, Cambridge, 1995.

\vspace*{-3mm}
\bibitem{Sha92}
J. O. Shallit.
Continued fractions with bounded partial quotients: A survey.
{\em Enseign. Math.}\ {\bf 38} (1992), 151--187.

\vspace*{-3mm}
\bibitem{Sol97}
B.~Solomyak.
Dynamics of self-similar tilings.
{\em Ergod.\ Th.\ \& Dynam.\ Syst.}\ {\bf 17}(1997), 695--738.

\vspace*{-3mm}
\bibitem{Sol98a}
B.~Solomyak.
Nonperiodicity implies unique composition for
self-similar translationally finite tilings.
{\em Discrete Comput.\ Geom.}\ {\bf 20} (1998), 265--279.

\vspace*{-3mm}
\bibitem{Sol98b}
B.~Solomyak.
Spectrum of dynamical systems arising from Delone sets.
{\em Quasicrystals and Discrete Geometry (Fields Institute Monographs, 10.)}
Ed.\ J.~Patera.
American Mathematical Society, Providence, RI 1998, pp.~265--275.

\vspace*{-3mm}
\bibitem{Thu89}
W.~Thurston. 
Groups, Tilings and Finite State Automata (AMS Colloquium Lecture Notes).
1989, unpublished.

\vspace*{-3mm}
\bibitem{Thu90}
W.~Thurston.
Conway's tiling groups.
{\em Amer.\ Math.\ Monthly} {\bf 97} (1990), 757--773.

\end{thebibliography}
\end{document}